\input amstex
\documentstyle{amsppt}
\magnification=\magstep1

\pageheight{9.0truein}
\pagewidth{6.5truein}

\NoBlackBoxes
\TagsAsMath

\loadbold

%%\long\def\ignore#1\endignore{\par DIAGRAM\par}
\long\def\ignore#1\endignore{#1}

\ignore
\input xy \xyoption{matrix} \xyoption{arrow}
          \xyoption{curve}  \xyoption{frame}
\def\edge{\ar@{-}}
\def\dttdar{\ar@{.>}}

\def\levelpool#1{\save [0,0]+(-3,3);[0,#1]+(3,-3)
                  **\frm<10pt>{.}\restore}
\def\dashedge{\ar@{--}}

\def\dropvert#1#2{\save+<0ex,#1ex>\drop{#2}\restore}
\endignore

\def\la{{\Lambda}}
\def\lamod{\Lambda\text{-}\roman{mod}}
\def\Lamod{\Lambda\text{-}\roman{Mod}}
\def \len{\operatorname{length}} 
\def\AA{{\Bbb A}}

\def\PP{{\Bbb P}}
\def\SS{{\Bbb S}}

\def\NN{{\Bbb N}}

\def\hom{\operatorname{Hom}}
\def\aut{\operatorname{Aut}}
\def\Aut{\operatorname{Aut}}

\def\stab{\operatorname{Stab}}

\def\top{\operatorname{top}}

\def\GL{\operatorname{GL}}

\def\K{\operatorname{K}}

\def\Im{\operatorname{Im}}

\def\autlap{\operatorname{Aut}_\la(P)}
\def\autlat{\operatorname{Aut}_\la(T)}

\def\End{\operatorname{End}}

\def\pr{\operatorname{pr}}

\def\B{{\Cal B}} 
\def\C{{\frak C}}
\def\BC{\bold{C}}

\def\bd{{\bold{d}}}

\def\F{{\Cal F}}

\def\I{{\Cal I}}

\def\K{{\Cal K}}
\def\T{{\Cal T}}
\def\L{{\Cal L}}

\def\m{{\frak m}}

\def\S{{\sigma}}
\def\s{{\frak s}}

\def\T{{\Cal T}}
\def\t{{\frak t}}

\def\zhat{\widehat{z}}

\def\GRASS{\operatorname{\text{\smc{grass}}}}

\def\Hom{\operatorname{Hom}}

\def\flag{\operatorname{\frak{Flag}}}
\def\Schu{\operatorname{\Schu}}

\def\Gr{\operatorname{Gr}}
\def\boldgrasstd{\operatorname{\frak{Grass}}^T_{\bold{d}}}
\def\maxtopdeg{\operatorname{\frak{Max}}^T_{\bold{d}}}
\def\imaxtopdeg{\operatorname{\frak{Max}}^{T_i}_{\bold{d_i}}}
\def\maxmoduli{\operatorname{\frak{ModuliMax}}^T_{\bold{d}}}
\def\imaxmoduli{\operatorname{\frak{ModuliMax}}^{T_i}_{\bold{d_i}}}
\def\maxmoduliM{\operatorname{\frak{ModuliMax}}(M,T')}
\def\maxmoduliMtop{\operatorname{\frak{ModuliMax}}(M,\top M)}
\def\autlap{\aut_\Lambda(P)}
\def\autlat{\aut_\Lambda(T)}

\def\unirad{\bigl(\autlap \bigr)_u}

\def\grassS{\operatorname{\frak{Grass}}(\S)}

\def\grassSS{\operatorname{\frak{Grass}}(\SS)}

\def\grass{\operatorname{\frak{Grass}}}
\def\biggrass{\GRASS}

\def\degeneration{\le_{\text{deg}}}
\def\underbardim{\operatorname{\underline{\vphantom{p}dim}}}
\def\id{\operatorname{id}}

\def\Schu{\operatorname{Schu}}

\def\smallstarblock{ \left[\smallmatrix *&\cdots&*\\ \vdots&&\vdots\\ *&\cdots&*
\endsmallmatrix\right] }

\def\vsubseteq{\hbox{$\bigcup$\kern0.1em\raise0.05ex\hbox{$\tsize|$}}}

\def\generic{{\bf 1}}
\def\codes{{\bf 2}}
\def\GeomII{{\bf 3}}
\def\Bor{{\bf 4}}
\def\Hil{{\bf 5}}
\def\domino{{\bf 6}}
\def\GeomI{{\bf 7}}
\def\menace{{\bf 8}}
\def\class{{\bf 9}}
\def\degen{{\bf 10}}
\def\hierarchies{{\bf 11}}
\def\Hum{{\bf 12}}
\def\King{{\bf 13}}
\def\New{{\bf 14}}
\def\Rie{{\bf 15}}
\def\Zwara{{\bf 16}}

\topmatter

\vskip0.3truein

\title Top-stable degenerations of finite dimensional representations II
\endtitle

\rightheadtext{Top-stable degenerations of finite dimensional representations II}

\author  H. Derksen, B. Huisgen-Zimmermann, J. Weyman
\endauthor

\address Derksen: Department of Mathematics, University of Michigan, Ann Arbor, MI 48109\endaddress

\email hderksen\@umich.edu \endemail

\address Huisgen-Zimmermann: Department of Mathematics, University of California, Santa Barbara, CA 93106 \endaddress

\email birge\@math.ucsb.edu\endemail

\address Weyman: Department of Mathematics, Northeastern University, Boston, MA 02115\endaddress

\email 	
j.weyman\@neu.edu\endemail

\thanks All three authors were partially supported by NSF grants while working on this project.  \endthanks

\abstract  Let $\la$ be a finite dimensional algebra over an algebraically closed field.  We exhibit slices of the representation theory of $\la$ that are always classifiable in stringent geometric terms.  Namely, we prove that, for any semisimple object $T \in \lamod$, the class of those $\la$-modules with fixed dimension vector (say $\bd$) and top $T$ which do not permit any proper top-stable degenerations possesses a fine moduli space.  This moduli space, $\maxmoduli$, is a projective variety.  Despite classifiability up to isomorphism, the targeted collections of modules are representation-theoretically rich:  indeed, any projective variety arises as $\maxmoduli$ for suitable choices of $\la$, $\bd$, and $T$.  In tandem, we give a structural characterization of the finite dimensional representations that have no proper top-stable degenerations. \endabstract

\endtopmatter

\thanks All three authors were partially supported by NSF grants.\endthanks

\document

%%%%%%%%%%%%%%%%%%%%%%%%%%%%%%
%%%%%%%%%%%%%%%%%%%%%%%%%%%%%%

\head 1. Introduction and preliminaries \endhead

 This article continues the study of degenerations begun in \cite{\degen}.
Let $\la$ be a finite dimensional algebra over an algebraically closed field $K$.  For convenience, we assume $\la$ to be basic, but our main results generalize to the non-basic situation.  Given any finite dimensional semisimple (left) $\la$-module $T$, we provide a structural description of those representations of  $\la$ which are maximal under the degeneration order among the representations $M$ with top $T$.  The top of $M$ is $M/JM$, where $J$ is the Jacobson radical of $\la$, and the degeneration order is the partial order on isomorphism classes defined by ``$M \degeneration N$ $\iff$  $M$ degenerates to $N$".  Moreover, we prove that the class $\C$ of $\la$-modules with fixed dimension vector $\bd$ which are degeneration-maximal among those with top $T$ possesses a fine moduli space, labeled $\maxmoduli$.  The points of $\maxmoduli$ encode (unique) normal forms to which the minimal projective presentations of the modules in $\C$ can be reduced.  In other words, our classification of $\C$ is in terms of normalized projective presentations, which form a family  satisfying a strong geometric universality condition; see ``families of modules" below. 

The representations we target are the simplest among those with fixed top and dimension vector, in the sense that their structure cannot be further unraveled along the degeneration order without enlarging the top.  However, despite classifiability in strict terms, their representation theory is complex.  This is attested to by the fact that arbitrary projective varieties arise as moduli spaces $\maxmoduli$ for suitable choices of $\la$, $\bd$, and $T$.  
\medskip

For somewhat more detail, we excerpt two readily stated portions of our main results, Theorems 3.5 and 4.4, and Corollary 4.5.  The dimension vector $\bd$ is kept fixed throughout the ensuing preview, and ``module" means ``finite dimensional left $\la$-module".  Moreover, we fix a full sequence $e_1, \dots, e_n$ of primitive idempotents of $\la$.  

The structural description of the modules which are degeneration-maximal among those with fixed top $T$ will arise as a consequence of the following straightforward geometric observation:  Namely, if $P$ is a projective cover of a module $M$ (equivalently, of the top $T$ of $M$), say $M \cong P/C$, then $M$ is devoid of proper top-stable degenerations if and only if the stabilizer of $C$ in $\aut_\la(P)$ is a parabolic subgroup of $\aut_\la(P)$.  In particular, the following result backs a rule of thumb that already surfaced in \cite{\degen}: namely, that non-existence of certain types of degenerations of $M$ is tantamount to invariance of $C$ under corresponding classes of endomorphisms of $P$.

\proclaim{Theorem A.  Structural characterization}  The module $M$ has no proper top-stable degenerations if and only if the following conditions are satisfied:

{\rm (a)}  $M \cong \bigoplus_{i=1}^n \bigoplus_{j=1}^{t_i} \la e_i/ C_{ij}$, where, for each $i \le n$, the $C_{ij}$ are left ideals contained in $\la e_i$ which are linearly ordered under inclusion.

{\rm (b)} $f(C_{ij}) \subseteq C_{kl}$ for all $f \in \Hom_\la(\la e_i, J e_k)$ and all pairs $(i,j)$ and $(k,l)$ {\rm {(not necessarily distinct)}}. \endproclaim

Thus, the modules without proper top-stable degenerations are direct sums of local modules (that is, of modules with simple tops) whose first syzygies satisfy the invariance condition spelled out under (b) of the theorem.  This invariance can be expressed more compactly in the following form, which is convenient for checks in concrete instances:  $\dim_K \Hom_\la(P, JM) = \dim_K \Hom_\la (M,JM)$. For reinforcements of the linkage between degenerations of $P/C$ and invariance properties of $C$ in $P$, see Observation 3.1, Theorem 3.6, and Corollary 3.8.

\proclaim{Theorem B.  Classification}  For any semisimple  $T \in \lamod$, the modules of dimension vector $\bd$ which are degeneration-maximal among those with top $T$ have a  fine moduli space, $\maxmoduli$,  that classifies them up to isomorphism.  The variety $\maxmoduli$ is projective. 

Moreover, given any module $M$ whose top $M/JM$ is contained in $T$, the closed subvariety of $\maxmoduli$ consisting of the points that correspond to degenerations of $M$ is a fine moduli space for the maximal top-$T$ degenerations of $M$. 
\endproclaim

As a consequence, on starting with any module $M$, one hits analyzable strata as one cuts through the full poset of degenerations of $M$ along successive enlargements of the top.
\medskip

How can the variety $\maxmoduli$ be realized in representation-theoretic terms? Here is an outline: Depending on a choice of Borel subgroup in the automorphism group $\autlat$ of the top $T$,  the variety $\maxmoduli$ can be located as a closed subset in the following projective variety $\boldgrasstd$ which parametrizes the isomorphism classes of modules with top $T$ and dimension vector $\bd$ (see \cite{\GeomII} for more detail): If we view the projective cover $P$ of $T$ as a $K$-vectorspace and write $|\bd|$ for the sum of the entries of $\bd$, the variety $\boldgrasstd$ is a closed subvariety of the classical Grassmannian $\Gr\bigl(\dim P - |\bd| ,\, P \bigr)$ of all $(\dim P - |\bd|)$-dimen\-sion\-al subspaces of $P$.  Namely, $\boldgrasstd$ consists of those points $C \in \Gr\bigl(\dim P - |\bd| ,\, P \bigr)$ which are $\la$-submodules of $JP$ satisfying $\underbardim P/C =\bd$.  Clearly, $\boldgrasstd$ carries a morphic $\autlap$-action, whose orbits are in one-to-one correspondence with the isomorphism classes of $\la$-modules with top $T$ and dimension vector $\bd$, via the assignment $\autlap.C \mapsto [P/C]$.  On the set $\maxtopdeg$ consisting of the points $C \in \boldgrasstd$ which correspond to representations that are degeneration-maximal among the modules with top $T$ (=  the set of points which are stabilized by a parabolic subgroup of $\autlap$),  the $\autlap$-action is reduced to an $\autlat$-action, but may still have orbits of arbitrarily high dimension.  However, given any Borel subgroup $\B$ of $\autlat$, the points in $\maxtopdeg$ whose $\autlat$-stabilizer contains $\B$ are in one-to-one correspondence with the isomorphism classes of the degeneration-maximal modules with top $T$.  It is the pared-down variety made up of these points that qualifies as a fine moduli space for our classification problem.  To the extent required by our present objective, we will recall the conceptual background under ``Notation, conventions, and supporting facts" at the end of this section.    

Finally, we provide a rough guide through the paper. In Section 2, we exploit the semidirect product decomposition of the acting group, $\autlap \cong \unirad \rtimes \autlat$, where $\unirad$ is the unipotent radical of $\autlap$.  Namely, we separately consider unipotent degenerations  and torus-type degenerations of a module $M \cong P/C$. These are degenerations which are either solely due to the action of $\unirad$ on $\boldgrasstd$, or else due to the action of the maximal tori in $\autlap$; the latter are seen to be representative for degenerations resulting from the action of the reductive group $\autlat$ (Corollary 3.9).  While there is major overlap between these two types of degenerations, and while there usually exist top-stable degenerations of $M$ that are neither unipotent nor of torus type, their separate study provides useful insight into the mechanisms by which the $\autlap$-action on $\boldgrasstd$ ``triggers" degenerations.  In particular, we will find existence of arbitrary proper top-stable degenerations to be detectable by means of these two types alone.  In other words, a module $M$ does not have any proper top-stable degenerations in case $M$ has neither proper unipotent nor proper torus-type top-stable degenerations.  

In Section 3, we characterize the following three classes of objects in $\lamod$:  those without proper unipotent top-stable degenerations, those without proper top-stable degenerations of torus-type, and finally those that have no proper top-stable degenerations altogether.  The corresponding results can be found in Observation 3.1, Theorem 3.6, and Theorem 3.5, respectively.  Section 4 deals with classification, and in Section 5, we give examples, applying a general construction principle (which, in other guises, has been used before) to back up our claim concerning the ubiquity of the varieties $\maxmoduli$.
Section 6 contains a proof of the main classification result, Theorem 4.4, and the construction of the pertinent universal family of representations.  
\bigskip

\noindent{\it Notation, conventions, and supporting facts:}

The Jacobson radical of $\la$ will be denoted by $J$, and $L+1$ will be the Loewy length of $\la$.  Given any $\la$-module $M$, we write $\underbardim M$ for the dimension vector of $M$, and we call $\SS(M) = \bigl( J^l M/J^{l+1} M \bigr)_{0 \le l \le L}$ the {\it radical layering\/} of $M$; note that the top $T$ of $M$ equals the first entry of this sequence.  The radical layering may be viewed as a numerical invariant of $M$, refining the dimension vector; it is given by the count of isomorphism classes of simple summands in the successive quotients $J^l M / J^{l+1} M$. 

We will systematically identify isomorphic semisimple modules.  A semisimple module $T$ is called {\it squarefree\/} if $T = \bigoplus_{i=1}^n S_i^{t_i}$ with $t_i \le 1$.  By a {\it semisimple sequence with dimension vector $\bd$\/}, we mean a sequence $\SS = (\SS_0, \SS_1, \dots, \SS_L)$ of semisimple $\la$-modules $\SS_l$, such that $\sum_{l=0}^L \underbardim \SS_l = \bd$;  we refer to the module $\SS_0$ as the {\it top\/} of the sequence $\SS$.  Moreover, we say that a module $M$ has {\it radical layering $\SS$\/} in case $\SS(M) = \SS$.  The following subset of $\boldgrasstd$ is a locally closed subvariety which parametrizes the isomorphism classes of modules with radical layering $\SS$ (with repetitions):
$$\grassSS = \{C \in \boldgrasstd \mid \SS(P/C) = \SS\},\  \  \  C \mapsto [P/C].$$
Clearly, the varieties $\grassSS$ are stable under the $\autlap$-action and partition $\boldgrasstd$.  For additional information about these varieties, see \cite{\class}, \cite{\generic}, and \cite{\hierarchies}.
\smallskip

\noindent{\bf Degenerations:} A {\it degeneration\/} of a module $M$ with dimension vector $\bd$ is any module $M'$ which is represented by a point in the closure of the $\GL(\bd)$-orbit of $M$ in the classical affine variety parametrizing the modules with dimension vector $\bd$.  For general background on degenerations, we refer to \cite{\Rie} and \cite{\Zwara}.  A degeneration $M'$ of $M$ is called {\it proper\/} if it is not isomorphic to $M$; it is called {\it top-stable\/} if $M'/JM' = M/JM$, {\it layer-stable\/} if $\SS(M') = \SS(M)$.  Suppose $C \in \boldgrasstd$.  From \cite{\degen} we recall that the top-stable degenerations of $M = P/C$ are precisely the modules $M' = P/C'$, where $C'$ traces the orbit closure $\overline{\autlap.C}$ in $\boldgrasstd$.  Analogously, the points $C'$ in the relative closure of $\autlap.C$ in $\grassSS$ yield the layer-stable degenerations $P/C'$ of $P/C$.   
\smallskip

\noindent{\bf Notation for path algebras modulo relations and their modules:}  Since $\la$ is a basic finite dimensional algebra
over an algebraically closed field,
we do not lose generality in assuming that $\la$  is a path
algebra modulo relations, say $\la = KQ/I$ for a quiver $Q$ and an admissible ideal $I \subseteq KQ$.   
The product $pq$ of two paths $p$ and $q$ in
$KQ$ stands for ``first $q$, then $p$"; in particular, $pq$ is zero
unless the end point of $q$ coincides with the starting point of $p$.  A {\it path in $\la$\/} is a residue
class of the form $p + I$, where
$p$ is a path in $K Q \setminus I$; we will suppress the residue
notation, provided there is no risk of ambiguity.  Further, we will
gloss over the distinction between the left $\la$-structure of a module
$M \in \Lamod$ and its induced $KQ$-module structure when there is no
danger of confusion.

The vertices $e_1, \dots, e_n$ of $Q$ will be
identified with the full sequence of primitive idempotents of $\la$ given by the corresponding paths of length zero in $Q$.  The representatives $\la e_i/ Je_i$ of the simples in $\lamod$ will be denoted by  $S_i $.  An element
$x$ of a module $M$ will be called a {\it top element\/} of $M$ if $x \notin JM$ and $x$ is {\it normed\/} by some $e_i$, meaning
that $x = e_i x$.  Any collection $x_1, \dots, x_t$ of top elements of $M$
generating
$M$ and linearly independent modulo $JM$ will be referred to as a {\it 
sequence of top elements of $M$\/}.
\smallskip

\noindent {\bf $\PP^1$-curves in $\boldgrasstd$\/}:  This installment of preliminaries is relevant only towards proofs and examples.   Let $C$ be any point in $\boldgrasstd$ and $\phi$ a morphism $X \rightarrow \autlap.C$, where $X$ is a dense open subset of $\AA^1$.  Projectivity of $\boldgrasstd$ guarantees the existence of a unique extension to a curve $\overline{\phi}: \PP^1 \rightarrow \overline{\autlap.C}$.  By $\lim_{\tau \rightarrow \infty} \phi(C)$ we denote the unique point in $\overline{\phi}(\PP^1 \setminus \AA^1)$.  Clearly, this latter point determines a top-stable degeneration of $P/C$.  Conversely, every top-stable degeneration of $P/C$ can be obtained in this fashion.  Indeed, in light of the fact that all orbit closures $ \overline{\autlap.C}$ are unirational projective varieties, this follows from a general curve-connectedness result of Koll\'ar (see \cite{\degen,  Proposition 3.6}).
\smallskip

\noindent {\bf Families of modules parametrized by varieties and moduli spaces\/}: 
The appropriate concept of a {\it family of $d$-dimensional $\la$-modules parametrized by a variety $X$\/} was introduced by King in \cite{\King}.  It is a pair $(\Delta, \delta)$, where $\Delta$ is a geometric vector bundle of rank $d$ over $X$ and $\delta: \la \rightarrow \End(\Delta)$ is a $K$-algebra homomorphism.  Our notion of {\it equivalence of families\/} parametrized by the same variety $X$ is finer than King's in general; it is the coarsest to separate isomorphism classes of modules:  namely, we set ``$(\Delta_1, \delta_1) \sim (\Delta_2, \delta_2)$" provided that, for each $x \in X$, the fibre of $\Delta_1$ over $x$ is $\la$-isomorphic to the fibre of $\Delta_2$ over $x$.  Given a family $(\Delta, \delta)$ parametrized by $X$ and a morphism of varieties $\tau: Y \rightarrow X$, we denote by $\tau^*(\Delta, \delta)$ the induced family parametrized by $Y$, that is, the pullback of $(\Delta, \delta)$ along $\tau$.  A family $(\Delta, \delta)$ of modules, from a class $\C$ say, is said to be {\it universal\/} (relative to the specified equivalence relation) in case any family of modules in $\C$ parametrized by some variety $Y$ is induced from $(\Delta, \delta)$ by a unique morphism $Y \rightarrow X$, where $X$ again denotes the parametrizing variety of the reference family.  Finally, $X$ is said to be a {\it fine moduli space\/} for the isomorphism classes of the modules in $\C$ if there exists a universal family parametrized by $X$.  
\smallskip

\noindent {\bf Non-standard conventions\/}: We fix a semisimple module $T \in \lamod$ with dimension vector $\t = (t_1, \dots ,t_n)$ and total dimension $t = | \t |  =  \sum_{i=1}^n t_i$, as well as a projective cover $P$ of $T$, together with a sequence $z_1, \dots,  z_t$ of top elements of $P$.  This latter sequence will be referred to as {\it the distinguished sequence of top elements of $P$\/}; the primitive  idempotent norming the top element $z_r$ will be denoted by $e(r)$.  

Unless explicitly revoked, we will identify $\autlat$ with a subgroup of $\autlap$, via the following {\it distinguished embedding\/} of $\autlat$ into $\autlap$.  Suppressing the residue notation $z_r + JP$ in $P/JP$, we will also view the $z_r$ as a $K$-basis for $T$;  that is, we view $T$ as the subspace $\bigoplus_{r=1}^t K z_r$ of $P$.  Accordingly, we identify an automorphism $g$ of $T$ with the automorphism of $P$ which sends $z_r$ to $g(z_r)$ for $r \le t$.  In other words, $\autlat$ will be regarded as the group of those automorphisms of $P$ which leave the subspace $\bigoplus_{r=1}^t K z_r$ invariant.  It is straightforward that $\autlap = \unirad \rtimes \autlat$ under this identification, as $\unirad = \{\id_P +  f \mid f \in \Hom_\la(P,JP)\}$.  Clearly, $\autlat \cong \prod_{i=1}^n \GL_{t_i}(K)$.
\smallskip

\noindent {\bf Top elements of $P$ versus tori in $\autlap$\/}: Since any two maximal tori in $\autlap$ are conjugate, these tori correspond bijectively to the equivalence classes of sequences of top elements of $P$, under the following equivalence: $(y_1, \dots, y_t) \sim (y'_1, \dots, y'_t)$ if  $y'_i \in K^* \tau(y_i)$ for all $i$, where $\tau$ is a permutation of $1,\dots,t$ such that $e(i)= e(\tau(i))$.  The correspondence is determined by the map which assigns to a sequence $(y_1, \dots, y_t )$ of top elements of $P$ the maximal torus $\T = \{ g \in \autlap \mid g(y_i)\in K^*y_i \text{\ for\ } i=1,\dots,t\}$ in $\autlap$.   The maximal torus determined by the distinguished sequence $(z_1, \dots, z_t)$ of top elements in this manner will be referred to as the {\it distinguished maximal torus of $\autlap$} and denoted by $\T_0$; note that $\T_0$ is a subgroup of $\autlat = \GL \bigl( \bigoplus_{r=1}^t K z_r \bigr)$ under the above identifications.   As it coincides with the group of diagonal matrices relative to the basis $(z_r)_{r \le t}$, we write its elements simply in the form $(a_1, \dots, a_t)$ with $a_r \in K^*$.  

For background concerning algebraic groups we refer to \cite{\Bor} and \cite{\Hum}, for the theory of moduli spaces to \cite{\New}.

%%%%%%%%%%%%%%%%%%%%%%%%%%%%%%
%%%%%%%%%%%%%%%%%%%%%%%%%%%%%%

\head 2.  Unipotent degenerations and torus-type degenerations 
\endhead

The following concepts were introduced in broader generality in  \cite{\hierarchies}.  Here, we trim them down to the top-stable scenario.

\definition{Definition 2.1}  Let
$P = \bigoplus_{r=1}^t \la z_r$ with $\la z_r \cong \la e(r)$ be
the distinguished $\la$-projective cover of $T$, and  $M \cong P/C$ with $C \in
\boldgrasstd$. Given a subset $X \subseteq \boldgrasstd$, we denote by
$\overline{X}$ its closure in $\boldgrasstd$.
\smallskip

\noindent $\bullet$ Given any subgroup $H \le \autlap$, we call a module
$M'$ a {\it top-stable degeneration of type $H$ of $M$\/} if $M' \cong P/C'$ for some point $C'$ in the orbit closure $\overline{g^{-1} Hg.C}$ of some conjugate
$g^{-1}Hg$ of $H$ in $\autlap$.  Moreover, such a degeneration $M'$ is
called a {\it proper} degeneration of type
$H$ of $M$ if $M' \not\cong M$, that is, if $C' \in \overline{g^{-1}
Hg.C}\setminus \autlap.C$ for some $g \in \autlap$.
\smallskip

Suppose that $M' \cong P/C'$ is a top-stable degeneration of type $H$ of $M$.  
 
$\bullet$  In the special case where $H = \unirad$, we also call $M'$
a  {\it unipotent top-stable degeneration of $M$\/}. (Note:  Since $\unirad$ is normal in $\autlap$, this means that $C' \in \overline{\unirad.C}$.)

$\bullet$  In the special case where $H$ is a torus, the distinguished maximal torus in $\autlap$, we also call $M'$ 
a {\it top-stable degeneration of $M$ of torus type\/}.      
\enddefinition

\definition{Remarks 2.2} {\bf A.} These definitions do not depend on the choice of the point $C$ in the $\autlap$-orbit  corresponding to $M$ in $\boldgrasstd$.   Indeed,  $M'$ is a top-stable degeneration of type $H$ of $M$ if and only if there exist points $D$ and $D'$ in the $\autlap$-orbits corresponding to $M$ and $M'$, respectively, such that $D' \in \overline{H.D}$.
Since every maximal torus in $\autlap$ is conjugate to the distinguished torus $\T_0$, this means in particular:  $M'$ is a torus-type top-stable degeneration of $M$ if and only if there are points $D$ and $D'$ in $\boldgrasstd$ such that $D' \in \overline{\T_0.D}$, $M \cong P/D$ and $M' \cong P/D'$. 
  
{\bf B.} For a subgroup $H$  of $\autlap$ and any point $C \in \boldgrasstd$, the following conditions are equivalent:
\smallskip

\noindent $\qquad \bullet$  $P/C$ has no proper top-stable degenerations of type $H$.
\smallskip
 
\noindent $\qquad \bullet$  $\overline{g^{-1}H g. C} \subseteq \autlap.C$ for all $g \in \autlap$.
\smallskip

\noindent $\qquad \bullet$ $ \overline{H.D} \subseteq \autlap.D$ for all points $D \in \autlap.C$.
\enddefinition

Roughly speaking, the unipotent degenerations of $M = P/C$ can be thought of as resulting from curves in $\autlap$ which ``pull composition factors of $M$ towards the top", while degenerations of torus type result from curves that apply ``torque"  to the structure of $M$ by placing increasingly discrepant weights on the direct summands of an indecomposable decomposition of $P$.  However, 
in general, the class of unipotent top-stable degenerations of $M$ has a large overlap with the torus-type degenerations of $M$.  We illustrate this fact with an easy example to which we will repeatedly refer in the sequel. 

\definition{Example 2.3}  Let $\la = KQ/\langle \omega^2 \rangle$, where $Q$ is
the quiver with two vertices, $1$ and $2$, and two arrows, a loop $\omega$
at the vertex $1$ and an arrow $\alpha$ from $1$ to $2$.   We
let $T = S_1^2$ and $P = \la z_1 \oplus \la z_2$, where the distinguished top elements
$z_1$ and $z_2$ are both normed by $e_1$.  Set $\bd = (3,2)$, and consider $M = P/C$ for $C = Jz_2$ $\in$
$\boldgrasstd$.   

If $C' = \la \alpha \omega
z_1 \oplus \la \omega z_2$, then $P/C'$ is a proper unipotent degeneration of $M$.   Indeed $C' =
\lim_{\tau \rightarrow \infty} g_{\tau}(C)$, where $g_\tau \in \unirad$
takes $z_1$ to $z_1$ and $z_2$ to $z_2 + \tau \omega z_1$ for $\tau \in \AA^1$.   

To see that $P/C'$ is also a degeneration of $M$ of torus type, we express $C'$ in the alternative form  $C' = \lim_{\tau \rightarrow \infty} f_\tau(C)$, where, for $\tau \in K^*$,  the automorphism $f_{\tau}$ of $P$ is determined by the assignments $z_1 \mapsto \tau z_1$ and $z_2 \mapsto z_2 + \omega z_1$.  Thus $C' \in
\overline{\T.C}$, where $\T$ is the torus subgroup of $\autlap$
determined by the sequence $(z_1, z_2 + \omega z_1)$ of top elements of
$P$. \enddefinition

For a first example of a top-stable degeneration which
is neither unipotent nor of torus type, we refer to
\cite{\degen, Example 5.8}. The top-stable degeneration
depicted in the central entry of the second column of Figure 4  --  call it $M'$  --  provides an
instance; in light of Proposition 2.4 and Corollary 2.7 below, this follows from the fact that $M'$ is a proper  layer-stable degeneration which is indecomposable.   In general, there is a plethora of such hybrid
degenerations, even in the case where $T$ is squarefree; consult \cite{\degen, Example 5.9}, again using the mentioned results. 

The proper unipotent degenerations $M'$ of a module $M$ always disturb the radical layering.  Namely, if $\SS = \SS(M)$ and $\SS' = \SS(M')$, then $\SS$ is properly dominated by $\SS'$ (written $\SS < \SS'$) in the following sense:  $\bigoplus_{l=0}^m \SS_l$ is a submodule of $\bigoplus_{l=0}^m \SS'_l $ for all $m \le L$, and inequality occurs for at least one $m$.  This was shown in  \cite{\hierarchies, Theorem 4.3}  in a more general context.  We record an immediate consequence which will be used in the next section.  

\proclaim{Proposition 2.4} Let $\SS$ be a semisimple sequence and $C \in \grassSS$.  Then the orbit $\unirad.C$ is closed in $\grassSS$.   In particular, proper unipotent degenerations fail to be layer-stable.
\endproclaim

All top-stable degenerations of a local module are necessarily
unipotent, since for a simple top $T$, the $\autlap$-action on $\boldgrasstd$
boils down to an $\unirad$-action.  So, in particular, all proper top-stable degenerations of a local module $M$ disturb the radical layering of $M$.

Proposition 2.4 tells us, in particular, that  $\overline{\unirad.C} \setminus \unirad.C$ has empty intersection with $\autlap.C$. There is no comparable result for degenerations of torus type, as is illustrated by the following example.  In fact, the situation where
$$\T.C\  \subsetneqq \  \overline{\T.C} \  \subseteq  \  \autlap.C$$
is a common occurrence.  

\definition{Example 2.5}
Let $\la = KQ$, where $Q$ is the quiver $1 \rightarrow 2$.  Moreover, let $T = S_1^2$, $\bd = (2,1)$, $P = \la z_1 \oplus \la z_2 \cong (\la e_1)^2$, and $C = Jz_2$.  Then clearly $\autlap.C$ is closed in $\boldgrasstd$; on the other hand, if $\T \subseteq \autlap$ is the torus subgroup  determined by the sequence $( z_1 + z_2, z_2)$ of top elements of $P$, we have $J z_1 \in \overline{\T.C} \setminus \T.C$. 
\enddefinition

We assemble a number of auxiliary facts concerning orbit structure and dimensions, which were established in \cite{\class, Proposition 2.9}.

\definition{Observation 2.6} Let $C \in
\boldgrasstd$.  As before, $t = \dim_K T$ is the number of simple summands of
$T$, meaning  that the maximal torus subgroups of
$\autlap$ are $t$-dimensional. 
\smallskip

{\rm (a)}  $\unirad.C \cong \AA^{\m}$, where 
$$\m = \dim_K \Hom_\la(P, JP/C) - \dim_K\Hom_\la(P/C, JP/C).$$
\smallskip

{\rm (b)}  Let $\T$ be a maximal torus subgroup
of $\autlap$.  As a variety, the orbit $\T.C$ is a torus of
dimension at most $t - 1$ and at least
$t - \s$, where
$\s$ is the number of indecomposable summands of $P/C$.  Moreover, there
exists a point $D \in \autlap.C$ such that $\T.D \cong (K^*)^{t -
\s}$.
\smallskip

{\rm (c)}  The orbit $\autlap.C$ has
dimension 
$$\dim_K  \Hom_\la (P,P/C) - \dim_K \End_\la(P/C).$$
\smallskip
\enddefinition

Both unipotent and torus-type degenerations are comparatively accessible in concrete instances.   In particular, the varieties $\overline{\unirad.C} \cong \AA^{\m}$ and $\overline{\T.C}$ for tori $\T$ in $\autlap$ are rational, and hence all points in their closures can be connected to $C$ by $\PP^1$-curves.

Part (b) of Observation 2.6 shows in particular that, for any
indecomposable module
$M = P/C$ and any maximal torus subgroup $\T \le \autlap$, the orbit $\T.C$
has dimension $t - 1$, where $t = \dim T$.  Consequently, all proper torus-type
degenerations of $M$ are decomposable.  Indeed, if $C' \in
\overline{\T.C} \setminus \T.C$, then $\dim \T.C' < \dim \T.C$, and  the inequality $\dim \T.C' \le t - 2$ guarantees decomposability of $P/C'$.  We record this for easy reference.

\proclaim{Corollary 2.7} Let $M \in \lamod$. Every proper top-stable degeneration of $M$ of torus type is decomposable. 
\qed \endproclaim  

 On the other hand, if $M$ is decomposable, proper top-stable degenerations of $M$ of torus type need not have a larger number of indecomposable summands than $M$.
 
\definition{Example 2.8}
Suppose $\la = KQ$, where $Q$ is the Kronecker quiver $\xymatrixcolsep{4pc} \xymatrix{
1 \ar@/^0.4pc/[r]^{\alpha} \ar@/_0.4pc/[r]_{\beta} &2 }$.   Let $T = S_1^2$, $\bd = (2, 2)$, $P = \la z_1 \oplus \la z_2 \cong (\la e_1)^2$, and $C = \la \alpha z_1 \oplus \la \beta z_2$ in $\boldgrasstd$.  We will check that $\la e_1 \oplus S_1$ is a top-stable degeneration of torus type of $M = P/C$:   Choose the torus $\T \le
\autlap$ which consists of the maps of the form $f_{(a_1, a_2)}$ for
$(a_1, a_2) \in (K^*)^2$, defined by
$z_1 + z_2 \mapsto a_1 (z_1 + z_2)$ and $z_2 \mapsto a_2 z_2$ with $a_i \in
K^*$.  Then $Jz_2 = \lim_{\tau \rightarrow \infty} f_{(1, \tau)}(C)$ belongs
to $\overline{\T.C}$. 
\enddefinition

The next proposition upgrades Corollary 2.7.

\proclaim{Proposition 2.9}  Let $M \in \lamod$, and let $U \subseteq M$ be
a top-stably embedded submodule, meaning that $JU = U \cap JM$. Then the top-stable degeneration
$U \oplus M/U$ of $M$  is of torus type.  

More precisely:  Suppose $\rho: P/C \rightarrow M$ is an isomorphism and $y_1, \dots, y_t$  a sequence of top elements of $P$ with the property that $\rho(y_1), \dots, \rho(y_s)$ is a sequence of top elements of $U$ for suitable $s \le t$.  If  $\T \le \autlap$ is the maximal torus subgroup corresponding to $y_1,\dots,y_t$ {\rm {(}}see Section 1, conventions{\rm {) }}, then  
$$U \oplus M/U \cong P/C' \ \ \text{for some} \ \  C' \in \overline{\T.C}.$$ 
 \endproclaim

\demo{Proof}  It is harmless to assume that $M = P/C$ with $C \in \boldgrasstd$.  Given top-embeddability of $U$ in $M$, we do not lose any generality
in assuming that $z_r + C$ for $r \le s$, form a full sequence of top elements
of $U$, where $s$  is the dimension of the top of $U$.  In other words, we may assume that the sequence of top elements $(y_r)_{r \le t}$ coincides with the distinguished sequence $(z_r)_{r \le t}$.  In particular, this adjustment entails $\T = \T_0$.  For
$\tau \in K^*$, consider the automorphism $\psi_{\tau} \in
\T_0$, defined by $z_i \mapsto z_i$ for $1 \le i \le s$ and
$z_i \mapsto \tau z_i$ for $s+1 \le i \le t$.   

We wish to prove that $C' : = \lim_{\tau \rightarrow \infty} \psi_\tau
(C)$ equals $(C \cap Q) \oplus \pi(C)$, where $Q = \bigoplus_{r=1}^s \la z_r$ and
$\pi: P \rightarrow Q' = \bigoplus_{r=s+1}^t \la z_r$ is the
canonical projection along $Q$; once this is established, we are done,
since, by construction, $C'$ belongs to the closure of $\T_0.C$ in
$\boldgrasstd$, $U \cong Q/ (C \cap Q)$, and $M/U \cong \pi(P)/\pi(C)$.  Let
$c_1,\dots, c_u$ be a $K$-basis for $C \cap Q$, and supplement it to a
$K$-basis $c_1, \dots, c_v$ of $C$, where
$u \le v$.  For each $k \le v$,
the limit $\lim_{\tau \rightarrow \infty} \psi_\tau (Kc_k)$ is a 
one-dimensional subspace of
$C'$.  It equals $Kc_k$ for $k
\le u$, since $\psi_\tau(c_k) = c_k$ in this case; in particular, $C
\cap Q \subseteq C'$. On the other hand, for $k \ge u+1$, we obtain
$\lim_{\tau \rightarrow \infty}
\psi_\tau(K c_k)$ $=$
$\lim_{\tau \rightarrow \infty} K \bigl[ 1/\tau \bigl((\id_P - \pi) (c_k) \bigr) + \pi(c_k) \bigr]$
$=$ $K \pi(c_k)$; here we invoke the computational rules developed in
\cite{\degen, Section 4.B}.   Since $\pi(C)$ is generated by the $\pi(c_k)$, we thus obtain the
claimed format of $C'$.        
\qed \enddemo

The argument backing Proposition 2.9 moreover yields the
following decomposition result.

\proclaim{Corollary 2.10}  Let $M \cong P/C$ with $C \in \boldgrasstd$, and let $\T$ be any maximal torus in $\autlap$.  Moreover, suppose that $y_1, \dots, y_t$ is the sequence of top elements of $P$ which corresponds to $\T$ in the sense of Section 1.  Then there exists a point
$C' \in \overline{\T.C}$ such that $C' = \bigoplus_{r=1}^t (C
\cap \la y_r)$.  In particular, the torus-type top-stable degeneration
$$M' = P/C' \cong \bigoplus_{r=1}^t \la y_r/ (C' \cap \la y_r)$$
of $M$ is a direct sum of local modules.
\endproclaim

%%%%%%%%%%%%%%%%%%%%%%%%%%%%%%
%%%%%%%%%%%%%%%%%%%%%%%%%%%%%%

\head 3.  
The structure of representations without proper top-stable degenerations  
\endhead

By definition, a module $M = P/C$ with $C \in \boldgrasstd$ is devoid of proper top-stable degenerations precisely when the orbit $\autlap.C$ is closed in $\boldgrasstd$; in other words, precisely when the stabilizer subgroup of $C$ in $\autlap$ is a parabolic subgroup of $\autlap$.
We will recoin this stabilizer condition in terms of a structural description of the corresponding module $M$.  For added information on the different mechanisms leading to degenerations, we separately focus on the cases where $M$ is either (A)
devoid of proper unipotent top-stable degenerations,  or (B) devoid of proper top-stable degenerations of torus type.  Each of the nonexistence conditions (A) and (B) is strictly weaker than nonexistence of arbitrary proper top-stable degenerations in general.  (Easy examples can be found below Observation 3.1 and Corollary 3.3.)  Yet, as we shall find, the combination of (A) and (B) is strong enough to preclude the existence of arbitrary proper top-stable degenerations of $M$.  

While a characterization of the modules satisfying (A) readily carries over from the special case where the top $T$ is squarefree   --  in this special case a characterization is already available  --   a useful description of the modules subject to (B) is harder to come by.  In light of Corollary 2.10, a necessary condition for the absence of proper torus-type top-stable degenerations is immediate however:  they are necessarily direct sums of local modules, that is $t - \s = 0$, where $t = \dim T$ and $\s$ is the number of indecomposable direct summands of $M$.

Recall that, in case the top $T$ of $M \cong P/C$ is squarefree, the isomorphism invariants
$\m$ and $t - \s$ of Observation 2.6  are known to govern the
size of the set of all top-stable degenerations of $M$ (see
\cite{\degen, Theorems 4.4 and 5.1}).  The first quantity,
$$\m = \dim_K \Hom_\la(P, J P/C) - \dim_K \Hom_\la(P/C, J P/C),$$
measures the deviation of $C$ from being invariant under homomorphisms
$P \rightarrow JP$, and hence the deviation of $C$ from being invariant under the maps in $\unirad$.  For a squarefree top $T$, vanishing of both invariants, $\m$ and $t - \s$, characterizes the situation where $M$ has no proper top-stable degenerations.  We will see that, on allowing multiple simple summands in $T$, these conditions will remain necessary for the absence of proper top-stable degenerations of $M$, but will no longer be sufficient.  

Clearly, the condition that $\unirad.C$ be a singleton precludes the existence of proper unipotent top-stable degenerations of $P/C$.  By the following observation, this condition is also necessary.  Typically, it is straightforward to check it in concrete instances by dint of condition (3) below.

\proclaim{Observation 3.1}  Let $M \cong P/C$ with $C \in \boldgrasstd$.  Then
the following statements are equivalent: 

\roster
\item"{\rm (1)}" $M$ has no proper unipotent top-stable degenerations.
\item"{\rm (2)}" $\unirad.C = \{C\}$, that is, $C$ is stable under all
maps in $\Hom_\la(P, JP)$.
\item"{\rm (3)}" $\m = 0$.
\endroster
\endproclaim

\demo{Proof }  The equivalence of  (2)  and  (3) follows from Observation 2.6(a), and the implication ``(2)$\implies$(1)" is trivial.  

Now asssume (1), meaning that the closure $\overline{\unirad.C}$ of $\unirad.C$ in $\boldgrasstd$ is contained in $\autlap.C$.   If $\SS = \SS(M)$, then clearly $\autlap.C \subseteq \grassSS$, and consequently $\overline{\unirad.C} \subseteq \grassSS$.  Since, by Proposition 2.4, $\unirad.C$ is closed in $\grassSS$, this forces $\unirad.C = \overline{\unirad.C}$ to be a projective variety.  On the other hand,  $\unirad.C$ is irreducible affine by Observation 2.6(a), whence we conclude that $\unirad.C$ is reduced to a point.  This proves (2).  \qed 
\enddemo

In Example 2.8, we described a module $M$ that satisfies the equivalent conditions of Observation 3.1, while having a proper top-stable degeneration of torus type.  In fact, any module of Loewy length $2$ is devoid of proper unipotent top-stable degenerations.

Next we analyze the situation where $M$ has no proper top-stable degenerations of torus type.  The outcome does not parallel that of the unipotent scenario:  The obviously sufficient condition that $\T.C$ be closed in $\boldgrasstd$ for all torus subgroups $\T$ of $\autlap$ clearly fails to be necessary in general (cf\. Example 2.5).   In a first installment, we will give a necessary condition for absence of proper torus-type top-stable degenerations of $M$ (Corollary 3.3).  This preliminary step will yield the main result of Section 3 (Theorem 3.5), which characterizes the modules of the section title.  Subsequently we will tie down the loose end concerning existence of degenerations of torus type.  

By a {\it flag variety \/} of a vector space $V$ we will always mean what is also labeled a {\it variety of partial flags\/}.  Our notation is as follows:  Given a sequence $1 < m_1 < m_2 < \cdots < m_u  = \dim V$, the corresponding flag variety $\flag(V, m_1, \dots, m_u)$ consists of the chains of subspaces $V_1 \subset \cdots \subset V_u$ of $V$ with $\dim V_i = m_i$.   Their relevance in the present context is due to the well-known fact that the quotients of $\autlat \cong \prod_{i=1}^n \GL_{t_i}$ by parabolic subgroups are direct products of flag varieties.

\proclaim{Lemma 3.2}  We view $\autlat$ as a subgroup of $\autlap$ by way of the distinguished sequence of top elements of $P$ (cf\. Section 1), and, given any subset $X \subseteq \boldgrasstd$, we denote by $\overline{X}$ its closure in $\boldgrasstd$.  Moreover, for $1 \le i \le n$, we let $\I_i$ be the set of all $r \in \{1, \dots, t\}$ such that $e(r) = e_i$.  

\noindent For any point $D \in \boldgrasstd$,  the following conditions {\rm {(i) --  (iii)}} are equivalent:
\medskip

{\rm (i)}   $\overline{\T.D} \subseteq \autlat.D$ for any torus subgroup $\T \subseteq \autlat$.
\smallskip

{\rm (ii)}  $\autlat.D$ is closed in $\boldgrasstd$, that is, the stabilizer $\stab_{\autlat}( D)$ is a parabolic subgroup of $\autlat$.
\smallskip

{\rm (iii)}  There exists a point $C \in \autlat.D$ which decomposes in the form
$$C = \bigoplus_{r=1}^t C \cap \la z_r,$$ 
such that this decomposition satisfies the following additional condition: If $C_r \le \la e(r)$ is the left ideal with $C \cap \la z_r = C_r z_r$, then, for each $i \in \{1, \dots, n\}$, the left ideals $\bigl(C_r\bigr)_{r \in \I_i}$ are linearly ordered under inclusion.
\medskip

\noindent  If {\rm {(i) -- (iii)}}  are satisfied, then $\autlat.D \cong \prod_{1 \le i \le n,\ t_i > 0} \F_i$, where the $\F_i$ are flag varieties depending only on the number of distinct left ideals in the family $(C_r)_{r \in \I_i}$ and the multiplicities of their occurrences.
\endproclaim

\demo{Proof}  
``(i)$\implies$(iii)".  Assume (i).  Then $\overline{\T_0.D} \subseteq \autlat.D$, where $\T_0$ is the distinguished torus in $\autlat$, and hence Corollary 2.10 yields a point $C \in \autlat.D$ such that $C = \bigoplus_{r=1}^t (C\cap \la z_r)$; say $C = g.D$ for a suitable $g \in \autlat$.   Moreover, we note that for any torus subgroup $\T$ of $\autlat$,
$$\overline{\T.C}  =  g\bigl( \overline{g^{-1} \T g \bigr).D} \subseteq \autlat.D = \autlat.C$$
by hypothesis, whence condition (i) carries over from $D$ to $C$.  We will check that this point $C$ satisfies the requirements listed in (iii).  Assume, to the contrary, that there exists  $i \in \{1, \dots, n\}$ together with indices $v,w \in \I_i$ such that the left ideals $C_v$ and $C_w$ in $\la e_i$ are not comparable.  To reach a contradiction, consider the maximal torus $\T_1$ in  $\autlat$ which is determined by the sequence of top elements $\zhat_1, \dots, \zhat_t$, where $\zhat_r = z_r$ for $r \ne v$ and $\zhat_v =  z_v + z_w$.  Using Proposition 2.9, one readily verifies that  $\overline{\T_1.C}$ contains a point $C'$ such that 
$$P/C' \cong \biggl( \bigoplus_{r \neq v,w} \la z_r/ C_r z_r \biggr)\  \oplus \ \la e_i/ (C_v \cap C_w) \ \oplus \la e_i / (C_v + C_w).$$  
Since $\la e_i/ (C_v \cap C_w) $ has higher dimension than either $\la e_i/C_v$ or $\la e_i/ C_w$, we conclude  $P/C' \not\cong P/C$.  This contradicts our hypothesis which postulates that $C' \in \autlat.C$. Thus the point $D$ satisfies condition (iii) .

``(iii)$\implies$(ii)".  Let $C \in \autlat.D$ be as in (iii).  Since the stabilizers of $C$ and $D$ in $\autlat$ are conjugate in $\autlat$, we only need to show that $\stab_{\autlat} C$ is a parabolic subgroup of $\autlat$.  It is therefore harmless to assume $C =D$.  We write $P_i = \bigoplus_{r \in \I_i} \la z_r$, and let $T_i = \bigoplus_{r \in \I_i} K z_r \cong (S_i)^{t_i}$ be the $S_i$-homogeneous component of $T$.  Moreover, $\bold{C}_i$ stands for  the direct sum $\bigoplus_{r \in \I_i} C_r z_r \subseteq P_i$.   Finally, we set $\bd_i = \underbardim P_i/\bold{C}_i$ and identify $\Aut_\la(T_i)$ with a subgroup in $\Aut_\la(P_i)$ as explained in Section 1.

In light of $C = \bigoplus_{i=1}^n \BC_i$, we have
$\stab_{\autlat}C  = \prod_{i=1}^n \stab_{ \Aut_\la(T_i)}\bold{C}_i$, and hence it suffices to prove that each of the factors
$\stab_{\Aut_\la (T_i)}\bold{C}_i$ is a parabolic subgroup of $\Aut_\la (T_i)$.  We will therefore not lose any generality in further simplifying our setup to the case where $T$ is homogeneous, say $T = T_1$, $P = P_1$ and $C = \bold{C}_1$.  Moreover, it is clearly harmless to reorder the
left ideals $C_r$ contained in $\la e_1$ so that      
$$C_1 = \cdots = C_{s_1} \supsetneqq C_{s_1 + 1} = \cdots = C_{s_1 + s_2}
\supsetneqq \cdots \supsetneqq C_{s_1 + s_2 + \cdots s_{u-1} + 1} = \cdots
= C_{s_1 + s_2 + \cdots s_u},$$
where $\sum_{j=1}^u s_j = t$.  In other words, the largest of
the left ideals $C_r \subseteq \la e_1$ occurs with multiplicity $s_1$,
etc., and the number of distinct left ideals among the $C_r$ is $u$.  Given that $T$ is homogeneous, our identification
of $\autlat$ with a subgroup of $\autlap$ amounts to an identification of $\autlat$ with $\GL(Kz_1 \oplus \cdots \oplus Kz_t)$.  Relative to the basis $(z_r)_{r \le t}$, the stabilizer
$\stab_{\autlat} C$ clearly equals the  subgroup 

$$\left[\matrix 
\smallstarblock &* &\cdots &\cdots &*\\
{\bold0} &\smallstarblock &* &\cdots &*\\
\vdots &&\ddots &\phantom{\smallstarblock} &\vdots\\
\vdots & &\phantom{\smallstarblock} &\smallstarblock &*\\
{\bold0} &\cdots &\cdots &{\bold0} &\smallstarblock
\endmatrix\right]$$

\noindent of $\autlat$, where the $j$-th block along the main diagonal belongs to $\GL_{s_j}(K)$.  Thus the subgroup  $\stab_{\autlat} C$ of $\autlat$ is indeed parabolic.  

Moreover, the orbit map $\autlat \rightarrow \autlat.C$ is separable, which guarantees that $\autlat.C \cong \autlat/\stab_{\autlat} C$.  Hence our supplementary claim becomes obvious on inspection of the stabilizer subgroup.
\smallskip

The implication ``(ii)$\implies $(i)" is trivial.  \qed
 \enddemo
 
We derive a necessary condition for absence of proper torus-type top-stable degenerations.  The corollary is actually a consequence of the proof of Lemma 3.2.  It hinges on the fact that the implication ``(i)$\implies$(iii)" in Lemma 3.2 remains valid if $\autlat$ is replaced by $\autlap$ wherever it occurs, the above argument following through, mutatis mutandis.

\proclaim{Corollary 3.3}  Suppose that $M$ has no proper top-stable degeneration of torus type, i.e., $M \cong P/D$ with $\overline{\T.D} \subseteq \autlap.D$ for all torus subgroups $\T$ of $\autlap$.  Then 
$$M = \bigoplus_{i=1}^n \bigoplus_{j=1}^ {t_i} M_{ij},$$
 where $M_{ij }\cong \la e_i/C_{ij}$, such that, for each $i \in \{1,\dots,n\}$, the left ideals $C_{ij} \subseteq \la e_i$ are linearly
ordered by inclusion.  
\endproclaim
 
\definition{Examples 3.4}
\noindent $\bullet$ The necessary condition  for absence of proper top-stable degenerations of torus type established in Corollary 3.3 fails to be sufficient in general.  Indeed, take $\la$, $T$, $\bd$, and $C \in \boldgrasstd$ as in Example 2.3.  Then $M = P/C$ is isomorphic to $\la e_1 \oplus S_1$ and thus clearly satisfies the conclusion of Corollary 3.3.  However, $M$ does have a proper degeneration of torus type, as we saw in Example 2.3. 
 \smallskip

\noindent $\bullet$  The absence of proper top-stable degenerations of torus type does not preclude the existence of proper unipotent top-stable degenerations.  Again, let $\la$ be the algebra of Example 2.3.  Then the local module $M = \la e_1/ \la \alpha$ has no proper torus type top-stable degenerations, but $M' = \la e_1 / \la \omega \alpha$ is a proper unipotent degeneration. 
\enddefinition 

We will exhibit a necessary and sufficient condition for the absence of proper torus-type degenerations at the end of the section.
At this point, we only need to tesselate the available pieces to obtain a
characterization of the modules of the section title.  
  
\proclaim{Theorem 3.5} Let $M$ be a module with dimension vector $\bd$ and top $T = \bigoplus_{i=1}^n S_i^{t_i}$.  Moreover, let $C$ be a point in $\boldgrasstd$ such that $M \cong P/C$.  
Then the following statements are equivalent:  

\roster
\item"{\rm (1)}" $M$ has no proper top-stable degeneration, i.e., the stabilizer subgroup $\stab_{\autlap} C$ is a parabolic subgroup of $\autlap$.
\item"{\rm (2)}" $M$ has no proper top-stable degeneration
which is either unipotent or of torus type.
\item"{\rm (3)}"  $M$ satisfies these two conditions:
\itemitem{$\bullet$}  $M$ is a direct sum of local modules, say $M =
\bigoplus_{i=1}^n \bigoplus_{j=1}^{t_i} M_{ij}$, where $M_{ij} \cong \la
e_i/C_{ij}$ {\rm (that is, $t -
\s = 0$)}, with the following additional property:
For each $i \le n$, the $C_{ij}$ are linearly ordered
under inclusion.
\itemitem{$\bullet$} $\dim_K \hom_\la(P, JM) - \dim_K \hom_\la(M,JM) = 0$.
\endroster

\noindent If conditions {\rm (1) -- (3)} are satisfied, then $\unirad$ stabilizes $C$, and $\autlap.C = \autlat.C$ is isomorphic to a direct product of flag
varieties $\F_i$,  where $\F_i$ depends only on the number of distinct left ideals in the family $(C_{ij})_{j \le t_i}$ and their multiplicities.
\endproclaim

\demo{Proof} The final assertion is covered by Lemma 3.2 and Corollary 3.3.  The implication ``(1)$\implies$(2)" is obvious.
 Towards verification of the remaining implications, we remark the following:  In light of the semidirect product decomposition $\autlap = \autlat \ltimes \unirad$, the condition $\unirad.C = \{C\}$ entails $\autlap.C = \autlat.C$.  Invoking Observation 3.1, we thus find that any of the conditions (1), (2), (3) guarantees equality of the orbits $\autlap.C$ and $\autlat.C$. Consequently,  the implications ``(2)$\implies$(3)$\implies$(1)" follow from Lemma 3.2.   \qed \enddemo

As a first application of Theorem 3.5, we supplement the necessary condition of Corollary 3.3 to a characterization of the modules without proper top-stable degenerations of torus type.   We already know that the first of the two conditions spelled out in Theorem 3.5(3) is necessary for absence of proper torus-type top-stable degenerations of $M$ (Corollary 3.3), whence we are looking for an appropriate relaxation of the second.

\definition{Remark}  Suppose that $M$ satisfies the first of the two conditions of Theorem 3.5(3).  Then the second is tantamount to the following invariance of the $C_{ij}$:  
\roster
\item"($\ddagger$)" \ \ If $(i,j)$ and $(k,l)$  are eligible pairs of indices, not necessarily distinct, and $f$ is a map in $\Hom_\la( \la e_i, J e_k)$, then $f(C_{ij}) \subseteq C_{kl}$. 
\endroster  \enddefinition

The correct weakening of $(\ddagger)$ for our present purpose is condition (2.b) in the upcoming theorem.  

\proclaim{Theorem 3.6}  Let $M$ be a module with top $T = \bigoplus_{i=1}^n S_i^{t_i}$ and dimension vector $\bd$.  Then the following conditions {\rm (1) -- (3)} are equivalent:
\smallskip

\roster
\item $M$ has no proper top-stable degeneration of torus type, that is,  $\overline{\T.C} \subseteq \autlap.C$ for all torus subgroups $\T$ of $\autlap$ and all points $C$ in the $\autlap$-orbit corresponding to $M$.
\smallskip

\item  $M \cong \bigoplus_{i=1}^n \bigoplus_{j=1}^{t_i} \la e_i / C_{ij}$, where the $C_{ij}$ are left ideals contained in $\la e_i$ such that 

\quad {\rm (2.a)} $C_{i1} \supseteq C_{i2} \supseteq \cdots \supseteq C_{it_i}$ for all $i$, 

 \noindent and

\quad {\rm (2.b)} $(i,j) \ne (k,l)$ implies $f(C_{ij}) \subseteq C_{kl}$ for all $f \in \Hom_\la(\la e_i, J e_k)$.
\smallskip

\item  $M \cong \bigoplus_{i=1}^n M_i$, where  each $M_i$ has top $S_i^{t_i}$, say $M_i \cong P_i/ \bold{C}_i$ with $P_i = \bigoplus_{z_r = e_i z_r} \la z_r$ and $\bold{C}_i \subseteq JP_i$, such that

\quad {\rm (3.a)} if $t_i \ge 2$, then $M_i$ does not have {\rm any} proper top-stable degeneration, and 

\quad{\rm (3.b)}  $i \ne k$ implies $f(\bold{C}_i) \subseteq  \bold{C}_k$ for all $f \in \Hom_\la(P_i, P_k)$.
\endroster
\endproclaim

We prepare for the proof of Theorem 3.6.  Again, we write $P = \bigoplus_{i=1}^n P_i$, where $P_i = \bigoplus_{r \le t,\, z_r = e_i z_r} \la z_r$.  In particular, the top of $P_i$ equals $S_i^{t_i}$. 

\proclaim{Lemma 3.7}   Suppose $C = \bigoplus_{i=1}^n \BC_i$ in $\boldgrasstd$, with $\BC_i \subseteq JP_i$ for $i \le n$, and let $M = P/C = \bigoplus_{i=1}^n P_i / \BC_i$.  Then the following conditions are equivalent: 
\smallskip

{\rm(i)}  Each point $D \in \autlap.C$ decomposes in the form $D = \bigoplus_{i=1}^n D \cap P_i$.
\smallskip

{\rm(ii)}   Each point $D \in \autlap.C$ decomposes in the form $D = \bigoplus_{i=1}^n f_i(\BC_i)$ for suitable $f_i \in \aut_\la(P_i)$.
\smallskip

{\rm(iii)}  $f(\BC_i) \subseteq \BC_k$ whenever $i \ne k$ and $f \in \Hom_\la(P_i, P_k)$.
\endproclaim

\demo{Proof}  The implication ``(i)$\implies$(iii)" follows from an application of (i) to the point $D = (\id_P + f).C$, and ``(ii)$\implies$(i)" is trivial.  To prove``(iii)$\implies$(ii)", suppose that (iii) holds, and let $D \in \autlap.C$.  Then $D = g (\id_P + h).C$, where $g \in \autlat$ and $h \in \Hom_\la (P, JP)$; this is due to the fact that $\autlap$ is the semidirect product of the groups $\autlat$ and $\unirad$, the former identified with a subgroup of $\autlap$ via the distinguished top elements $z_r$ of $P$ as usual.  By the definition of the embedding $\autlat \subseteq \autlap$, the map $g$ leaves the $P_i$ invariant, and therefore it suffices to check that $( \id_P + h).C$ has a decomposition as postulated in (ii).   We write the map $h$ in terms of its components, $h =  \sum_{i , k=1}^n h_{ik}$, where $h_{ik} \in \Hom_\la(P_i, JP_k)$.

First we ascertain that we do not lose any generality in assuming $h_{ii} = 0$ for all $i$.  Indeed $\id_P + h = (\id_P + \psi)(\id_P + \chi)$, where $\psi = \sum_{i=1}^n h_{ii}$ and 
$$\chi = \sum_{i,k, \, i \ne k} h_{ik}\  - \  \sum_{i,k, \, i \ne k} h_{kk}h_{ik}\ + \  \sum_{i,k, \, i \ne k} (h_{kk})^2 h_{ik} \  - + \ \cdots \ + (-1)^{L-1} \sum_{i,k, \, i \ne k} (h_{kk})^{L-1}  h_{ik};$$
given that $J^{L+1} = 0$, this follows from the fact that $h_{uv}(P_u)  \subseteq JP_v)$ for all $u,v$.
By construction, $\id_P + \psi$ leaves all the summands $P_i$ of $P$ invariant, whence it is enough to verify that $(\id_P + \chi).C$ has the desired decomposition property.  Since $\chi$ has no nontrivial components $\chi_{ii} \in \Hom_\la(P_i, JP_i)$, this justifies reduction to the case $h_{ii} = 0$. 

In this situation, condition (iii) guarantees $(\id_P + h)(C) = C$.  Thus $D$ decomposes as required. 
\qed \enddemo

\proclaim{Corollary 3.8}  Retain the hypotheses and notation of Lemma 3.7, and set $M_i = P_i /\BC_i$.  Then the following conditions are equivalent for $M = \bigoplus_{i=1}^n M_i$:
\smallskip

{\rm(i)}  Every top-stable degeneration $M'$ of $M$ is of the form $M' = \bigoplus_{i=1}^n M_i'$, where each $M_i'$ is a degeneration of $M_i$.
\smallskip

{\rm(ii)}  $f(\BC_i) \subseteq \BC_k$ whenever $i \ne k$ and $f \in \Hom_\la(P_i,P_k)$.
\smallskip

More specifically, every torus-type top-stable degeneration of $M$ is a direct sum of torus-type degenerations of the  $M_i$ in case {\rm(ii)} is satisfied.
\endproclaim

\demo{Proof}  The implication ``(ii)$\implies$(i)" is immediate from Lemma 3.7.  For the converse, suppose (ii) to fail; say there exists $f \in \Hom_\la(P_1, P_2)$ and $c \in \BC_1$ such that $f(c) \notin \BC_2$.  For $\tau \in \AA^1 \setminus \{0\}$, we define $f_{\tau} \in \Aut_\la(P)$ via $f_{\tau}(z_1) = \frac{1}{\tau}\, z_1$ and $f_{\tau}(z_i) = z_i$ for $i \ge 2$; moreover, we set $g_{\tau} = f_\tau \circ (\id_P + f)$.  Clearly, $\tau \mapsto g_{\tau}$ is in turn a morphism $\AA^1 \setminus \{0\} \rightarrow \autlap$, since $f(P_1) \subseteq J P_2$.  Setting $C = \bigoplus_{i=1}^n \BC_i$, we observe that $C' = \lim_{\tau \rightarrow \infty} g_\tau(C)$ contains $\bigoplus_{i=2}^n \BC_i$, as well as the one-dimensional space 
$$\lim_{\tau \rightarrow \infty} K g_\tau(c) = \lim_{\tau \rightarrow \infty} K \bigl(\frac{1}{\tau}\, z_1 + f(c) \bigr) = K f(c).$$
This shows that $\dim (P_2 \cap C') > \dim (P_2 \cap C)$, whence the degeneration $M' = P/C'$ of $M$ fails to satisfy (i).   

In light of Remarks 2.2, the supplementary statement is straightforward.  \qed 
\enddemo 

Contrasting Corollary 3.8, in general degenerations of direct sums need not even satisfy cancellation; see \cite{\Rie}.  
For use in the upcoming argument, we note that the proof of the corollary shows in particular: Any module $M = \bigoplus_{i} M_i$ with top($M_i$) homogeneous of type $S_i$, which violates condition (ii), does have a proper top-stable degeneration of torus type; indeed, in the notation of ``(i)$\implies$(ii)", the module $P/D$ with $D = (\id_P + f).C$ is isomorphic to $M$, and $C' = \lim_{\tau \rightarrow \infty} f_\tau(D)$.

\demo{Proof of Theorem 3.6}  (1)$\implies$(2).  We adopt (1).  Corollary 3.3 then tells us that $M$ decomposes in the form $M \cong \bigoplus_{i=1}^n \bigoplus_{j=1}^{t_i} \la e_i / C_{ij}$, where the $C_{ij}$ are left ideals contained in $\la e_i$ for $i \le n$, such that condition (2.a) is satisfied.  So we only need to verify condition (2.b).  As we just pointed out, the proof of Corollary 3.8 guarantees $f(C_{ij}) \subseteq C_{kl}$ whenever $i \ne k$ and $f \in \Hom_\la(\la e_i, \la e_k)$.  So it suffices to consider maps $f \in \Hom_\la (\la e_i, J e_i)$.  In the upcoming indirect argument for (2), we may thus assume that $M$ has  homogeneous top $S_i^{t}$.  Accordingly, we simplify our notation to $M = P/C$ with $C = \bigoplus_{r=1}^t C_r z_r$, where $C_r  = C_{ir} \subseteq \la e_i$. By (2.a), we have $C_1 \supseteq C_2 \supseteq \cdots \supseteq C_t$. Suppose there exist indices $v \ne w$ in $\{1, \dots, t\}$ and $f \in \hom_\la( \la e_i, Je_i)$ with $f(C_v) \not\subseteq C_w$. Then, clearly, $t \ge 2$ and $f(C_1) \not\subseteq C_t$.  From this setup we will construct a proper top-stable degeneration of $M$ of torus type.  To this end, we abbreviate $f(e_i) \in e_i J e_i$ to $\omega$.  Let $\T$ be the torus subgroup of $\autlap$ which corresponds to the sequence $(\zhat_1, \dots , \zhat_t)$ of top elements of $P$, where $\zhat_1 = z_1 + \omega z_t$  and $\zhat_r = z_r$ for $2 \le r \le t$.  Defining $g_\tau \in \T$ for $\tau \in \AA^1 \setminus \{0\}$ via $g_\tau(\zhat_r) = \zhat_r$ for $r < t$ and $g_\tau (\zhat_t) = \tau \zhat_t$, we obtain $C' = \lim_{\tau \rightarrow \infty} g_\tau (C) \in \overline{\T.C}$.  Thus $P/C'$ is a torus-type degeneration of $M$.

To ascertain that $P/C' \not\cong M$, we also realize $P/C'$ as a unipotent degeneration $P/C''$ of $M$.  For this purpose, we consider the curve $\varphi: \AA^1 \rightarrow \unirad.C$,  $\, \tau \mapsto h_\tau.C$, where $h_\tau \in \unirad$ sends $z_1$ to $z_1 - \tau \omega z_t$  and $z_r$ to $z_r$ for $r > 1$.   The extension of $\varphi$ to a curve $\PP^1 \rightarrow \overline{\unirad.C}$ is in turn denoted by $\varphi$.  
 We observe that $\varphi$ is not constant.  Indeed, let $c \in C_1$ with $f(c) \notin C_t$; then $\varphi(0) = C$, while $\varphi(1) \ne C$ since $h_1(c z_1) =  c z_1 -  f(c) z_t \notin C$ by construction.  Consequently, the irreducible projective variety $\varphi(\PP^1)$ is not contained in the affine variety $\unirad.C$.  In light of $\varphi(\AA^1) \subseteq \unirad.C$, this yields $C'' = \lim_{\tau \rightarrow \infty} h_\tau.C \in \overline{\unirad.C} \setminus \unirad.C$.  Finally, we check that $C'' = C'$:   Indeed,  
 $$g_\tau(z_1) = g_\tau(\zhat_1 - \omega z_t) = \zhat_1 - \tau \omega z_t = z_1 - (\tau - 1) \omega z_t = h_{\tau - 1}(z_1),$$
and $g_\tau(z_r) = h_{\tau - 1}(z_r)$ for $2 \le r < t$, whereas $g_\tau(z_t) =  \tau h_{\tau - 1}(z_t)$. Since $\tau \ne 0$, this yields $g_\tau(C_r z_r) = h_{\tau - 1}(z_r C_r)$ for all $r \in \{1, \dots, t\}$, and thus $g_\tau(C) = h_{\tau - 1}(C)$, showing $C' = C''$. In light of Proposition 2.4, the degeneration $P/C'$ hence fails to be layer-stable; in particular, $P/C' \not\cong M$.  This  
contradiction completes the proof of ``(1)$\implies$(2)". 

(2)$\implies$(3).  Assuming (2), we set $\BC_i =  \bigoplus_{r \le t,\, z_r = e_i z_r} C_r z_r$ and $M_i = P_i/\BC_i$.  Condition (3.b) is immediate from (2.b).  To verify (3.a), suppose $t_i \ge 2$.  By Theorem 3.5, it suffices to supplement the inclusion $f(C_{ij}) \subseteq C_{ij}$ for all $j \in \{1, \dots t_i\}$ and any map $f \in \Hom_\la ( \la e_i, Je_i)$.  But, since $f(C_{i1}) \subseteq C_{i\, t_i}$ by (2.b), this is obvious in light of $C_{i1} \supseteq \cdots \supseteq C_{i\, t_i}$.  

(3)$\implies$(1).  Assume (3).  In light of Corollary 3.8, one gleans from (3.b) that any top-stable degeneration of M of torus type is a direct sum of torus-type top-stable degenerations of the $M_i$.  Hence it suffices to show that the individual $M_i$ do not have any proper torus-type top-stable degenerations. If $t_i \ge 2$, then this is a consequence of (3.a), and if $t_i = 0$, there is nothing to show.  For $t_i = 1$, the module $M_i$ is local, and thus the tori in $\Aut_\la(P_i)$, copies of $K^*$, act trivially on the Grassmannians of subspaces of $P_i$.  Therefore, $M_i$ is devoid of proper torus-type top-stable degenerations  in this case as well.   \qed
  \enddemo 
  
The proof of ``(3)$\implies$(1)" of Theorem 3.6 in fact shows that the modules without proper top-stable degenerations of torus type satisfy the upcoming reinforced condition.  It further clarifies how maps in $\autlap$  trigger degenerations of $M$, by separating contributions from the reductive and unipotent parts of $\autlap$.  
 
\proclaim{Corollary 3.9}  Again suppose that $M \in \lamod$ has top $T$.  Then $M$ has no proper top-stable degenerations of torus type if and only if $M$ has no proper top-stable degenerations of type $\autlat$.  \qed \endproclaim
  
 While for $\dim T = 1$, all modules with top $T$ are devoid of proper top-stable degenerations of torus type, non-existence of such degenerations turns into a strong condition as the top grows. The case where all simple summands of $T$ occur with a multiplicity $> 1$ is addressed by the following combination of Theorem 3.6 and Lemma 3.7.
 
\proclaim{Corollary 3.10}  Suppose $T = \bigoplus_{i=1}^n S_i^{t_i}$ such that all nonzero multiplicities $t_i$ are at least $2$, and let $M$ be a module with top $T$.   Then absence of proper torus-type top-stable degenerations of $M$ implies absence of arbitrary proper top-stable degenerations.  \qed \endproclaim

Finally, we point out that condition (3.b) in Theorem 3.6 is not redundant, even if $t_i \ge 2$ for all $i \le n$.

\definition{Example 3.11}
Let $Q$ be the quiver $\xymatrixcolsep{4pc} \xymatrix{
1 \ar@/^0.4pc/[r]^{\alpha} &2  \ar@/^0.4pc/[l]^{\beta} }$. Moreover, let $\la = KQ/I$, where $I$ is the ideal generated by all paths of length $3$, and $M = S_1^2 \oplus (\la e_2)^2$.  Then $M_1 = S_1^2$ and $M_2 = (\la e_2)^2$ are both without proper top-stable degenerations.  Nonetheless, $M$ has proper top-stable degenerations of torus type; for instance, $M' = (\la e_1/ \la \beta \alpha)^2 \oplus (\la e_2/ \la \alpha \beta)^2$ is such a degeneration. 
\enddefinition

%%%%%%%%%%%%%%%%%%%%%%%%%%%%%%
%%%%%%%%%%%%%%%%%%%%%%%%%%%%%%

\head 4. Classification \endhead

In this section, we use Theorem 3.5 to identify a fine moduli space, $\maxmoduli$, for the representations which are degeneration-maximal among those with fixed dimension vector $\bd$ and top $T$.   We will find this moduli space to be a closed subvariety of $\boldgrasstd$.  Concrete incarnations of $\maxmoduli$ within $\boldgrasstd$ depend on a choice of Borel subgroup of $\autlat$.  

More precisely, we will locate $\maxmoduli$ inside the following subvariety of $\boldgrasstd$, namely, 
$$\maxtopdeg \ := \ \{C \in \boldgrasstd \mid P/C \text{\ has no proper degeneration with top\ } T\}.$$  
In view of Theorem 3.5, $\maxtopdeg$ can alternatively be described as the set of all points $C \in \boldgrasstd$ with the following two properties:  $C$ is a fixed point under the action of $\unirad$, and $\stab_{\autlat} C$ is a parabolic subgroup of $\autlat$.  
In particular,  this implies that $\unirad$ acts trivially on $\maxtopdeg$, that is, the $\autlap$-action on $\maxtopdeg$ is reduced to an $\autlat$-action.  However, in general this action on $\maxtopdeg$ will still have large orbits.  Our task consists of showing that, in factoring the $\autlat$-action out of $\maxtopdeg$, we obtain a sufficiently well-behaved quotient.   Again, we will identify $\autlat$ with a subgroup of $\autlap$ by way of the distinguished top elements $z_r$ of $P$  (see Section 1). We start by recording that $\maxtopdeg$ is in turn projective, being a morphic image of the projective variety $G/ \B \times \{C \in \boldgrasstd | B \subseteq \stab_{\autlat} C\}$, where $B$ is any Borel subgroup of $\autlat$.

\proclaim{Observation 4.1}  $\maxtopdeg$ is a closed subvariety of $\boldgrasstd$. \qed \endproclaim

{\bf Our choice of Borel subgroup $\B$ of $\autlat$.}  Any Borel subgroup of $\autlat$ will serve our purposes.  But to ascertain unequivocal definitions (that is, in order to eliminate the uncertainty factor ``up to isomorphism induced by conjugation in $\autlat$"), we will fix a Borel subgroup of $\autlat$, based on a specific ordering of the distinguished top elements of $P$.  As before, $\t = (t_1, \dots, t_n)$ is the dimension vector of $T$.  Assume that the distinguished  sequence $z_1, \dots, z_t$ of top elements of $P$ is ordered in such a fashion that the first $t_1$ of the $z_r$ are normed by the idempotent $e_1$, the next $t_2$ are normed by $e_2$, etc.  Thus, setting $r_0 = 0$ and $r_i = t_1 + \cdots + t_i$ for $1 \le i \le n$, and defining
$$P_i =  \bigoplus_{r=r_{i-1} + 1}^{r_i} \la z_r,$$
we obtain $P = P_1 \oplus \cdots \oplus P_n$ with $P_i/JP_i \cong S_i^{t_i}$.   Accordingly,  $\autlat$ may be identified with $ \prod_{i=1}^n \GL(V_i)$, where 
$$V_i = \bigoplus_{r=r_{i-1} + 1}^{r_i} Kz_r.$$
This setup allows us to specify a Borel subgroup $\B$ of $\autlat$ as follows: 
$$\B = B_1 \times \cdots \times B_n,$$
where $B_i \subseteq \GL(V_i)$ is the group of upper triangular automorphisms relative to the basis $(z_r)_{r_{i-1} + 1 \le\, r \le\, r_i}$ of $V_i$.   

Accordingly, $\B \ltimes \unirad$ will be referred to as {\it the distinguished Borel subgroup of $\autlap$\/}.

\definition{Definition 4.2 and immediate observations} 
$$\maxmoduli   := \{ C \in \boldgrasstd  \mid \,  \B \ltimes \unirad \subseteq \stab_{\autlap} ( C)\}.$$

\noindent  In other words, $\maxmoduli$ is the set of all fixed points of the action of the distinguished Borel subgroup of $\autlap$ on $\boldgrasstd$.  In particular, $\maxmoduli$ is a closed (i.e., a projective) subvariety of $\boldgrasstd$.  Moreover,  $\maxmoduli$ is clearly contained in $\maxtopdeg$.  In fact, under our indexing conventions, 
$$\maxmoduli = \{ C = \bigoplus_{r=1}^t C_r z_r \in \maxtopdeg \mid C_{r_{i-1} + 1} \supseteq  \cdots \supseteq  C_{r_i} \ \text{for all}\ i \le n\}.$$
\enddefinition

The next observation shows $\maxmoduli$ to be a plausible candidate for a moduli space for those representations with top $T$ and dimension vector $\bd$ which have no proper top-stable degenerations.  Namely, the assignment ``$C \mapsto [P/C]$" is a one-to-one correspondence from the points of $\maxmoduli$ on one hand to the isomorphism classes of the specified modules on the other.  Indeed:

\proclaim{Observation 4.3} For each point $C \in \maxtopdeg$, there exists precisely one point $C' \in \maxmoduli$ such that $P/C \cong P/C'$. \endproclaim

\demo{Proof} Existence of $C'$ is due to the fact that $\stab_{\autlap}(C)$ is a parabolic subgroup of $\autlap$, whence $\stab_{\autlat}(C)$ contains a conjugate of $\B$ in $\autlat$.  For uniqueness, suppose $C', C'' \in \maxmoduli$ with $P/C \cong P/C'$.  This means that  $g.C' = C''$ for some $g \in \autlat$ and both $\stab_{\autlat}(C')$ and $\stab_{\autlat}(C'') = g\stab_{\autlat}(C') g^{-1}$  contain $\B$.  Therefore the groups $\stab_{\autlat}(C')$ and $ g\stab_{\autlat}(C') g^{-1}$ coincide (see, e.g., \cite{Borel, 11.17(i)}).  We conclude $g \in \stab_{\autlat}(C')$, because the stabilizer $\stab_{\autlat}(C')$  coincides with its normalizer in $\autlat$.  Hence $C'' = g.C' = C'$ as claimed.  \qed \enddemo

We thus obtain a map 
$$\pi: \maxtopdeg \longrightarrow \maxmoduli, \quad \quad C \mapsto C',$$
where $C'$ is chosen according to Observation 4.3.  By construction, $\pi$ is a surjection, the fibers of which coincide with the $\autlap$-orbits (= $\autlat$-orbits) of $\maxtopdeg$.  An explicit description of $\pi$  is as follows:  $\pi(C)$ is the unique point $C'$ in the $\autlat$-orbit of $C$ with the following properties: 
\roster
\item"$\bullet$" $C' = \bigoplus_{r=1}^t C' \cap \la z_r$, whence in particular, $C' \cap P_i  = \bigoplus_{j=r_{i-1} + 1}^{r_i}  C'_{ij} \, z_j$ for certain (uniquely determined) left ideals $C'_{ij} \subseteq \la e_i$, and
\item"$\bullet$" for $i \le n$, these left ideals satisfy 
$C'_{ij} \supseteq C'_{ik}$  whenever $j \le k$.
\endroster

\proclaim{Theorem 4.4}  The projective variety $\maxmoduli$ is a  fine moduli space for the isomorphism classes of representations which are degeneration-maximal among those with dimension vector $\bd$ and top $T$.

Supplementary information:  $\maxmoduli$ is the geometric quotient of the variety $\maxtopdeg$ modulo its $\autlap$-action \rm{(}which equals its $\autlat$-action{\rm)}, i.e.:

$\bullet$ $\pi: \maxtopdeg \rightarrow \maxmoduli$ is a surjective open morphism, whose fibers coincide with the $\autlap$-orbits of $\maxtopdeg$. 

$\bullet$  For every open subset $U$ of $\maxmoduli$, the comorphism $\pi^0$ of $\pi$ induces a $K$-algebra isomorphism from the ring $\Cal O (U)$ of regular functions on $U$ to the ring of those regular functions in $\Cal O \bigl( \pi^{-1} (U) \bigr)$ which are constant on the $\autlap$-orbits of $\pi^{-1} (U)$. 
\endproclaim

In non-technical terms, the corresponding universal family of $\la$-modules is the family $\bigl( P/C \bigr)$, where $C$ traces $\maxmoduli$.  By construction, it hits the isomorphism type of any module which is degeneration-maximal in our comparison class precisely once.   The universal property of this family is rooted in the following universal property of $\pi$, which is an immediate consequence of the fact that this map is a geometric quotient of $\maxtopdeg$ by $\autlap$:  Namely, given any morphism $\chi: \maxtopdeg \rightarrow Z$ which is constant on the $\autlap$-orbits, there exists a unique morphism $\chi': \maxmoduli \rightarrow Z$ with $\chi = \chi' \circ \pi$.  

We defer the proof of Theorem 4.4 to Section 6.  

\proclaim{Corollary 4.5}  Let $M$ be a module with top $T$, and $T'$ a semisimple module containing $T$.  Then the degenerations of $M$ which are maximal among those with top $T'$ possess a fine moduli space.  This moduli space, $\maxmoduliM$, is in turn projective.

In particular, the maximal top-stable degenerations of $M$ have a fine moduli space. \endproclaim

\demo{Proof of Corollary 4.5}  Let $\underbardim M = \bd$, and suppose that  $\underbardim T \le  \bd$ (if this inequality fails, our claim is void).  $P$ continues to denote the projective cover of $T$.  The variety in which we will locate a moduli space for the maximal top-$T$ degenerations of $M$ is $[\biggrass_\bd(\la)]_P$, as introduced in \cite{\hierarchies. Section 2}; it consists of the submodules $D$ of $P$ with ${\underline \dim}\, P/D = \bd$.  Note that we do not require $D \subseteq JP$, whence this variety parametrizes the isomorphism types of all $\la$-modules with dimension vector $\bd$ and top contained in $T$. In particular, $[\biggrass_\bd(\la)]_P$ includes $\boldgrasstd$ and contains a point $C$ such that $M \cong P/C$.  The intersection $\maxmoduli \cap \overline{\autlap.C}$ is the postulated moduli space. \qed \enddemo

\definition{Questions 4.6}
Theorem 4.4 and Corollary 4.5 prompt numerous questions, such as:  
\smallskip

{\bf(1)}  Which projective varieties arise as moduli spaces $\maxmoduli$ for suitable $\la$, $T$ and $\bd$?  

{\bf(2)} Which occur as $\maxmoduliM$ for prescribed $M$ and $T'$?  

{\bf(3)} Given $\la$, $T$ and $\bd$, how many irreducible components does $\maxmoduli$ have, how do their geometric characteristics depend on algebraic properties of $\la$, and how can their generic modules be accessed? 
\smallskip

 As announced earlier, the answer to the first question is ``all do"  (Example 5.4 below).  The answers to the other questions are open.  In Section 5, we analyze an example in Loewy length $3$ from the listed viewpoints (Example 5.1).  For this low Loewy length, we expect the example to be prototypical in the following sense:  namely, we expect all irreducible components of the varieties $\maxmoduliMtop$, for $M \in \lamod$, to be direct products of projective spaces; the components of the spaces $\maxmoduli$ are already far more complex, even in this example.  
 \enddefinition

%%%%%%%%%%%%%%%%%%%%%%%%%%%%%%
%%%%%%%%%%%%%%%%%%%%%%%%%%%%%%

\head 5. Examples \endhead

Our first  example illustrates the case $J^3 = 0$.  In this case, determining the spaces $\maxmoduliM$ for specific choices of $M$ from quiver and relations is a manageable task, purely combinatorial. We include some detail in our specific instance, instead of giving a procedural manual.  On the side we note: Extrapolating from the construction below, we find in particular that arbitrarily high numbers of irreducible components of $\maxmoduliM$ can be realized for $J^3 = 0$ and $\dim M/JM = \dim T = 2$.

\definition{Example 5.1} Let $\la = KQ/I$, where $Q$ is the quiver given below and $I \subseteq KQ$ is the ideal generated by all paths of length $3$, together with $\omega_i \omega_j$ for $1 \le i,j \le 6$, $\beta \omega_i$ for $1 \le i \le 3$, and $\alpha \omega_i$ for $4 \le i \le 6$.  

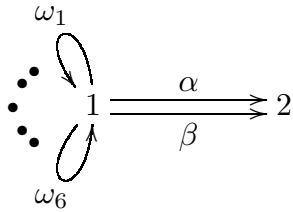
\ignore
$$\xymatrixcolsep{1.8pc}\xymatrixrowsep{2pc}
\xymatrix{
&  \\
 &&1 \ar@(u,ul)_{\omega_1} \ar@(dl,d)_{\omega_6}  
\ar@{}[ul]_(0.45){}="om1" \ar@{}[dl]^(0.45){}="om2"
\ar@{{}{*}}@/_1pc/"om1";"om2"
\ar[rr]<0.5ex>^{\alpha}
\ar[rr]<-0.5ex>_{\beta} &&2  \\
 &
}$$
\endignore

\noindent We let $T = S_1^2$ and $\bd = (12, 10)$. Moreover, we consider the left $\la$-module $M = P/C$ with top $T$, where $P = \la z_1 \oplus \la z_2$ with $\la z_1$ $\cong$ $\la z_2$ $\cong$ $ \la e_1$, and $C = C_1 z_1 \oplus C_2 z_2 \subseteq JP$ for the following left ideals $C_i \subseteq \la e_1$:  We take $C_1 = 0$, and, setting 
$$L = \la(\alpha \omega_1 + \alpha \omega_2) + \la (\beta \omega_4 + \beta \omega_5) + \la \omega_3 + \la \omega_6,$$ 
we define $C_2 = L + \la \alpha + \la \beta$.  Thus $M \cong \la e_1 \oplus (\la e_1/ C_2)$.  Optically more easily digestible is the following graph of $M$ (see \cite{\domino}  or \cite{\menace} for our graphing conventions):

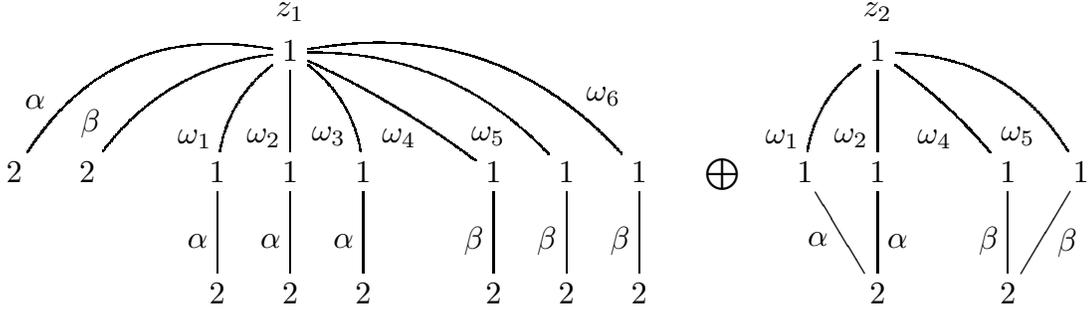
\ignore
$$\xymatrixcolsep{1.0pc}\xymatrixrowsep{2pc}
\xymatrix{
 &&&&&1 \dropvert3{z_1}
\edge@/_1.5pc/[1,-4]_(0.8){\alpha}
\edge@/_0.8pc/[1,-3]_(0.85){\beta} \edge@/_/[dl]_(0.8){\omega_1}
\edge[d]_(0.7){\omega_2} \edge@/^/[dr]_(0.7){\omega_3}
\edge@/^0.2pc/[1,3]_(0.6){\omega_4}
\edge@/^1pc/[1,4]_(0.75){\omega_5}
\edge@/^1.5pc/[1,5]^(0.8){\omega_6}  &&&&& && &1 \dropvert3{z_2}
\edge@/_/[dl]_(0.8){\omega_1} \edge[d]_(0.7){\omega_2}
\edge@/^0.2pc/[drr]_(0.6){\omega_4}
\edge@/^1pc/[1,3]_(0.75){\omega_5}  \\  
 &2 &2 &&1 \edge[d]_(0.55){\alpha} &1
\edge[d]_(0.55){\alpha} &1
\edge[d]_(0.55){\alpha} &&1 \edge[d]_(0.55){\beta} &1
\edge[d]_(0.55){\beta} &1 \edge[d]_(0.55){\beta} &\bigoplus &1
\edge[dr]_(0.45){\alpha} &1 \edge[d]^(0.55){\alpha} &&1
\edge[d]_(0.55){\beta} &1 \edge[dl]^(0.45){\beta}  \\
 &&&&2 &2 &2 &&2 &2 &2  &&&2 &&2
}$$
\endignore

The fine moduli space for the isomorphism classes of maximal top-stable degenerations of $M$, namely $\maxmoduliM$, has four irreducible components ${\bold{Comp}}_1, \dots, {\bold{Comp}}_4$, isomorphic to $\PP^1 \times \PP^1$, $\PP^1$, $\PP^1$, and $\PP^0$, respectively.   

We provide a formal description of the maximal top-stable degenerations of $M$ represented by the points in 
${\bold{Comp}}_1$. Given a point in $\PP^1 \times \PP^1$, written in the form 
$$\bigl([k_1:k_2:k_3], [k_4: k_5:k_6] \bigr) \ \  \text{with\ } k_1 = k_2 \text{\  and \ } k_4 = k_5,$$ 
the corresponding module is $\bigl(\la z_1/ D_1 z_1\bigr) \oplus \bigl( \la z_2/ D_2 z_2 \bigr)$, where $D_1$ is the $2$-dimensional left ideal $K(\sum_{i=1}^3 k_i \alpha \omega_i) + K(\sum_{i=4}^6 k_i \beta \omega_i)$ $\subseteq \la e_1$, and $D _2$ is the left ideal $L$ defined above.
Clearly, $D_1 \subseteq D_2$ and $f(D_i) = 0$ for $i = 1,2$ and all $f \in \Hom_\la(\la e_1, J e_1)$.  So the modules $P/D$ do not have any proper top-stable degenerations by Theorem 3.5.  We verify that they are degenerations of $M$:  For $\bigl([k_1:k_2: k_3] ,[k_4: k_5: k_6]\bigr)$ as above and $\tau \in K$, we let $g_\tau \in \autlap$ be defined by $z_1 \mapsto z_1$ and $z_2 \mapsto z_2 + \tau \sum_{i=1}^6  k_i \omega_i z_1$.  Then 
$$g_\tau (C) = L z_2 + K\bigl( \alpha z_2 + \tau \sum_{i=1}^3  k_i \alpha \omega_i z_1 \bigr) + K\bigl( \beta z_2 + \tau \sum_{i=4}^6  k_i \beta \omega_i z_1 \bigr),$$
whence $\lim_{\tau \rightarrow \infty} g_\tau(C)$ equals $D$.  For our notation, see Section 1, for a full justification of this computation, see \cite{\degen, Section 4.B}.

Explicit descriptions of the other irreducible components of $\maxmoduliM$ are obtained analogously.   To prove that all maximal degenerations of $M$ sharing the top $T$ with $M$ belong to $\bigcup_{1 \le i \le 4} {\bold{Comp}}_i$, one first argues that  all such degenerations have radical layering $\SS = (S_1^2, \, S_1^{10} \oplus S_2^ 4, \, S_2^6)$ and then uses Koll\'ar's curve connectedness (see ``$\PP^1$-curves in $\boldgrasstd$" in Section 1) to exclude the modules in $\grassSS$ that fail to be degenerations of $M$.  

Below, we display a graph of the general module corresponding to each of the irreducible varieties ${\bold{Comp}} _i$, to provide a more intuitive understanding of the universal family of modules parametrized by $\maxmoduliM$.
 In the graph pertaining to ${\bold{Comp}} _1$, the left-hand ``pool", outlined in dots, indicates that the space generated by the elements $\alpha \omega_i z_1$ for $i = 1,2,3$ has dimension $2$, a  linear dependence relation being $(k \alpha \omega_1 + k \alpha \omega_2 + k' \alpha \omega_3) z_1 = 0$ with $k,k' \in K^*$; the relation is communicated in terms of a point in $\PP^2$ adjacent to the pool.  More generally, a pool including $m$ vertices in the graph of a module $X$ signals that the corresponding elements in $X$ span a space of dimension $m-1$, while any choice of $m-1$ of these elements are linearly independent.   Note that, in general,  such a hypergraph  encodes the isomorphism class of the module $X$ only if every pool of $m \ge 2$ vertices is accompanied by a point in $\PP^m$.  The pool of the two elements $\alpha \omega_1 z_2$ and $\alpha \omega_2 z_2$ in the graph of the general module in ${\bold{Comp}} _1$ is determined by the relation $\alpha \omega_1 z_2 + \alpha \omega_2 z_2 = 0$ in $M$.  Analogously, $\beta \omega_4 z_2 + \beta \omega_5 z_2 = 0$ by the definition of $M$.

\ignore
$$\xymatrixcolsep{0.9pc}\xymatrixrowsep{2pc}
\xymatrix{
\dropvert3{{{\bold{Comp}}}_1:} &&&&1 \dropvert3{z_1}
\edge@/_1.5pc/[1,-4]_(0.8){\alpha}
\edge@/_0.8pc/[1,-3]_(0.85){\beta} \edge@/_/[dl]_(0.8){\omega_1}
\edge[d]_(0.7){\omega_2} \edge@/^/[dr]_(0.7){\omega_3}
\edge@/^0.2pc/[1,3]_(0.6){\omega_4}
\edge@/^1pc/[1,4]_(0.75){\omega_5}
\edge@/^1.5pc/[1,5]^(0.8){\omega_6}  &&&&& && &&&&1 \dropvert3{z_2}
\edge@/_1.5pc/[1,-4]_(0.8){\alpha}
\edge@/_0.8pc/[1,-3]_(0.85){\beta} \edge@/_/[dl]_(0.8){\omega_1}
\edge[d]_(0.7){\omega_2}
\edge@/^0.2pc/[drr]_(0.6){\omega_4}
\edge@/^1pc/[1,3]_(0.75){\omega_5}  \\  
2 &2 &&1 \edge[d]_(0.55){\alpha} &1
\edge[d]_(0.55){\alpha} &1
\edge[d]_(0.55){\alpha} &&1 \edge[d]_(0.55){\beta} &1
\edge[d]_(0.55){\beta} &1 \edge[d]_(0.55){\beta} &\bigoplus &2 &2
&&1 \edge[dr]_(0.45){\alpha} &1 \edge[d]^(0.55){\alpha} &&1
\edge[d]_(0.55){\beta} &1 \edge[dl]^(0.45){\beta}  \\
 &&&2 \levelpool2 &2 \dropvert{-5}{[k:k:k']} &2 &&2 \levelpool2 &2
\dropvert{-5}{[l:l:l']} &2  &&&&&&2 &&2 
}$$

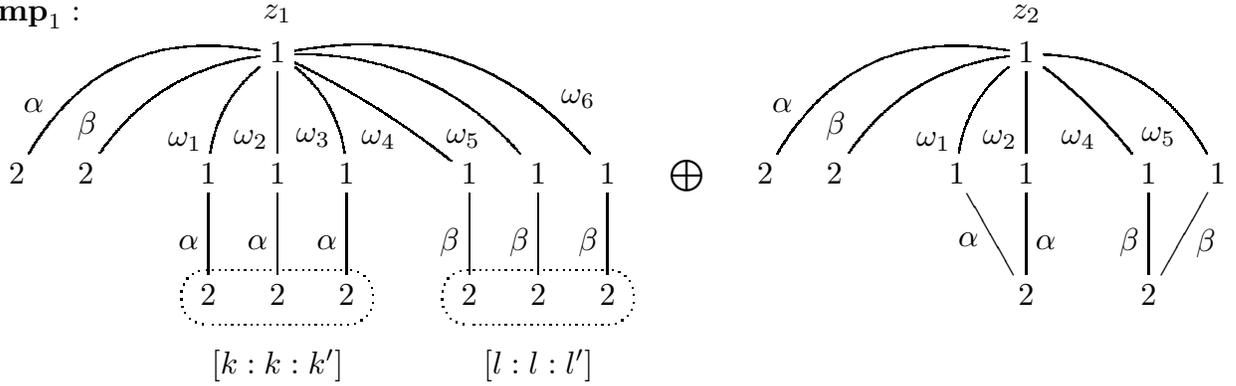
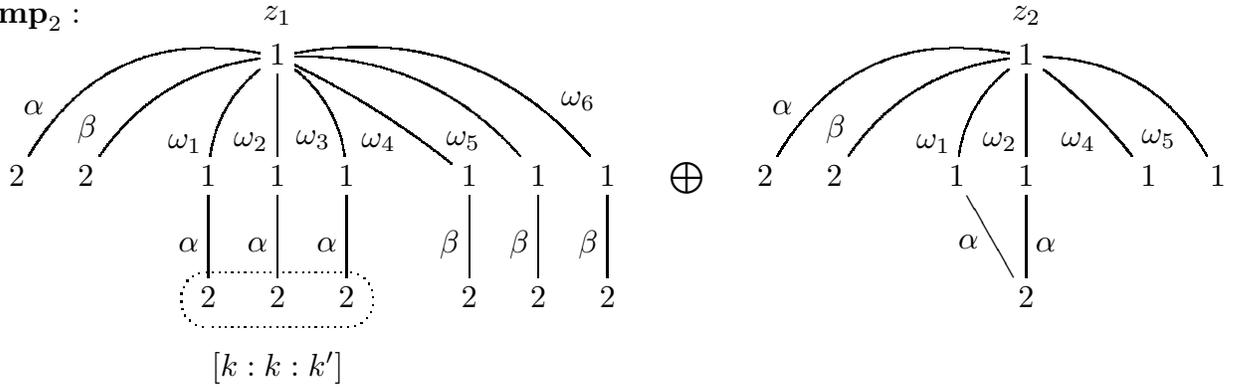
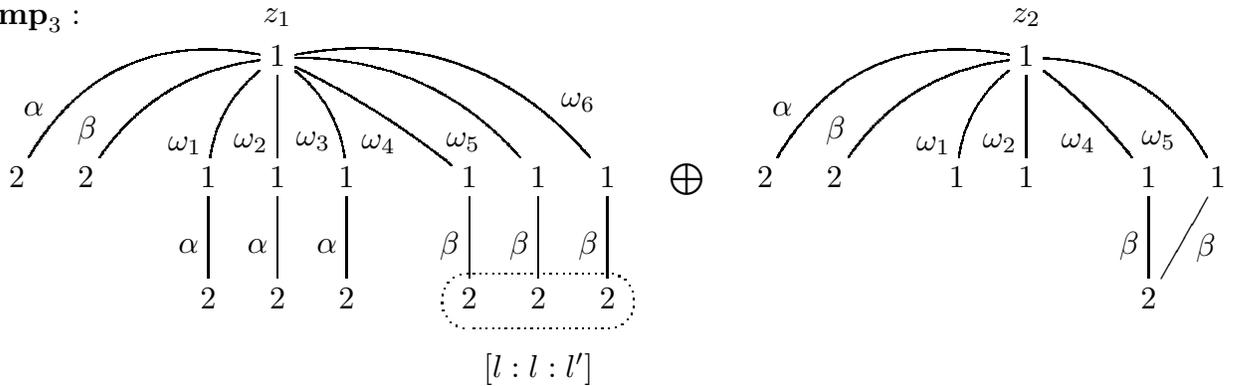
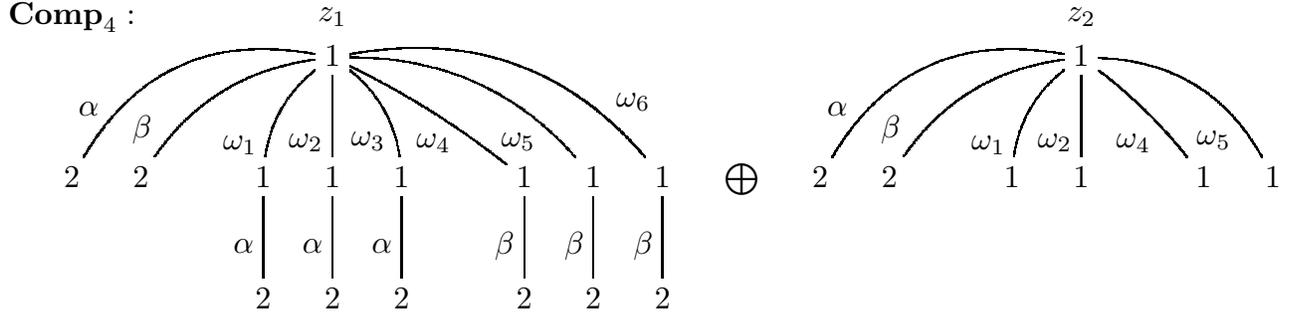
$$\xymatrixcolsep{0.9pc}\xymatrixrowsep{2pc}
\xymatrix{
\dropvert3{{{\bold{Comp}}}_2:} &&&&1 \dropvert3{z_1}
\edge@/_1.5pc/[1,-4]_(0.8){\alpha}
\edge@/_0.8pc/[1,-3]_(0.85){\beta} \edge@/_/[dl]_(0.8){\omega_1}
\edge[d]_(0.7){\omega_2} \edge@/^/[dr]_(0.7){\omega_3}
\edge@/^0.2pc/[1,3]_(0.6){\omega_4}
\edge@/^1pc/[1,4]_(0.75){\omega_5}
\edge@/^1.5pc/[1,5]^(0.8){\omega_6}  &&&&& && &&&&1 \dropvert3{z_2}
\edge@/_1.5pc/[1,-4]_(0.8){\alpha}
\edge@/_0.8pc/[1,-3]_(0.85){\beta} \edge@/_/[dl]_(0.8){\omega_1}
\edge[d]_(0.7){\omega_2}
\edge@/^0.2pc/[drr]_(0.6){\omega_4}
\edge@/^1pc/[1,3]_(0.75){\omega_5}  \\  
2 &2 &&1 \edge[d]_(0.55){\alpha} &1
\edge[d]_(0.55){\alpha} &1
\edge[d]_(0.55){\alpha} &&1 \edge[d]_(0.55){\beta} &1
\edge[d]_(0.55){\beta} &1 \edge[d]_(0.55){\beta} &\bigoplus &2 &2
&&1 \edge[dr]_(0.45){\alpha} &1 \edge[d]^(0.55){\alpha} &&1 &1  \\
 &&&2 \levelpool2 &2 \dropvert{-5}{[k:k:k']} &2 &&2 &2 &2  &&&&&&2
&& 
}$$

$$\xymatrixcolsep{0.9pc}\xymatrixrowsep{2pc}
\xymatrix{
\dropvert3{{{\bold{Comp}}}_3:} &&&&1 \dropvert3{z_1}
\edge@/_1.5pc/[1,-4]_(0.8){\alpha}
\edge@/_0.8pc/[1,-3]_(0.85){\beta} \edge@/_/[dl]_(0.8){\omega_1}
\edge[d]_(0.7){\omega_2} \edge@/^/[dr]_(0.7){\omega_3}
\edge@/^0.2pc/[1,3]_(0.6){\omega_4}
\edge@/^1pc/[1,4]_(0.75){\omega_5}
\edge@/^1.5pc/[1,5]^(0.8){\omega_6}  &&&&& && &&&&1 \dropvert3{z_2}
\edge@/_1.5pc/[1,-4]_(0.8){\alpha}
\edge@/_0.8pc/[1,-3]_(0.85){\beta} \edge@/_/[dl]_(0.8){\omega_1}
\edge[d]_(0.7){\omega_2}
\edge@/^0.2pc/[drr]_(0.6){\omega_4}
\edge@/^1pc/[1,3]_(0.75){\omega_5}  \\  
2 &2 &&1 \edge[d]_(0.55){\alpha} &1 \edge[d]_(0.55){\alpha} &1
\edge[d]_(0.55){\alpha} &&1 \edge[d]_(0.55){\beta} &1
\edge[d]_(0.55){\beta} &1 \edge[d]_(0.55){\beta} &\bigoplus &2 &2 &&1 &1 &&1 \edge[d]_(0.55){\beta} &1 \edge[dl]^(0.45){\beta} \\
 &&&2 &2 &2 &&2 \levelpool2 &2
\dropvert{-5}{[l:l:l']} &2  &&&&&& &&2 
}$$

$$\xymatrixcolsep{0.9pc}\xymatrixrowsep{2pc}
\xymatrix{
\dropvert3{{{\bold{Comp}}}_4:} &&&&1 \dropvert3{z_1}
\edge@/_1.5pc/[1,-4]_(0.8){\alpha}
\edge@/_0.8pc/[1,-3]_(0.85){\beta} \edge@/_/[dl]_(0.8){\omega_1}
\edge[d]_(0.7){\omega_2} \edge@/^/[dr]_(0.7){\omega_3}
\edge@/^0.2pc/[1,3]_(0.6){\omega_4}
\edge@/^1pc/[1,4]_(0.75){\omega_5}
\edge@/^1.5pc/[1,5]^(0.8){\omega_6}  &&&&& && &&&&1 \dropvert3{z_2}
\edge@/_1.5pc/[1,-4]_(0.8){\alpha}
\edge@/_0.8pc/[1,-3]_(0.85){\beta} \edge@/_/[dl]_(0.8){\omega_1}
\edge[d]_(0.7){\omega_2}
\edge@/^0.2pc/[drr]_(0.6){\omega_4}
\edge@/^1pc/[1,3]_(0.75){\omega_5}  \\  
2 &2 &&1 \edge[d]_(0.55){\alpha} &1
\edge[d]_(0.55){\alpha} &1
\edge[d]_(0.55){\alpha} &&1 \edge[d]_(0.55){\beta} &1
\edge[d]_(0.55){\beta} &1 \edge[d]_(0.55){\beta} &\bigoplus &2 &2
&&1 &1 &&1 &1  \\
 &&&2 &2 &2 &&2 &2 &2 
}$$
\endignore
\smallskip

We present only one of the irreducible components of the moduli space $\maxmoduli$ classifying {\it all\/} degener\-a\-tion-maximal modules with top $T = S_1^2$ and dimension vector $\bd = (12,10)$:  Namely, the unique irreducible component  of $\maxmoduli$ which includes both of the components ${\bold{Comp}}_1$ and ${\bold{Comp}}_2$ of $\maxmoduliM$. We denote it by ${\bold{Comp}}_0$.

To describe it, we abbreviate the subspace $\sum_{i=1}^6 K \omega_i$ of $\la e_1$ by $X$, and the subspace $\sum_{i=1}^3 K \alpha \omega_i +  \sum_{i=4}^6 K \beta \omega_i$ by $Y$.  The component ${\bold{Comp}}_0$ is isomorphic to the subvariety of  
$\Gr(1,Y) \times  \Gr(2,X)  \times \Gr(5,Y)$ which consists of the points  $(U,V,W)$ satisfying $U + \alpha V + \beta V \subseteq W$.  The module corresponding to such a triple is $\bigl( \la z_1 / U z_1\bigr) \oplus \bigl( \la z_2/ (V + W) z_2 \bigr)$.  \qed
\enddefinition

The following class of examples results from \cite{\class, Corollary 4.5}, but can also be deduced from King's work in \cite{\King}, as is detailed in \cite{\class, ahead of Corollary 4.5}.  

\proclaim{Example 5.2}  Let $\la = KQ/I$, where $Q$ is a quiver without oriented cycles. 

For a simple top $T$, say $T = S_1$, and any dimension vector $\bd$, we then have
$$\maxmoduli = \boldgrasstd. \ \ \ \ \ \qed$$
\endproclaim

The conclusion of Example 5.2 remains true under the following weakened blanket hypothesis:  For every oriented cycle $c$  in the quiver $Q$, the power $c^{\len(c)+1}$ belongs to $I$.   On the side we remark that, to some extent, $\maxmoduli$ is amenable to computation.  Namely, polynomial equations for the charts of a ``representation-theoretically distinguished" affine cover of $\boldgrasstd$ can be obtained algorithmically from quiver and relations of $\la$; the algorithm has been implemented in \cite{\codes}.

We consider a somewhat larger supply of examples akin to the preceding one. 

\proclaim{Example 5.3}  Again suppose that $\la = KQ/I$, where $Q$ has no oriented cycles, but now assume that $T$ is a homogeneous semisimple module, say of type $S_1$ and dimension $t$.  Then $\maxmoduli$ is the union of the following intersections, with $\bd_1 \le \bd_2 \le \cdots \le \bd_t$ tracing the sequences of dimension vectors with $\sum_{r=1}^t \bd_r = \bd$:
$$\bigl( \grass^{S_1}_{\bd_1} \times \cdots \times \grass^{S_1}_{\bd_t} \bigr) \ \cap \ \flag\bigl({\la e_1, m_1, \dots, m_t}\bigr),$$  
where $m_r = \dim_K \la e_1 - |\bd_r|$. \qed
\endproclaim

\demo{Proof}  From Example 5.2, we know that the moduli space of the degeneration-maximal modules with top $S_1$ and any dimension vector $\bd_r$ equals $\grass^{S_1}_{\bd_r}$.  Moreover, our hypothesis on $\la$ ensures that $\Hom_\la (\la e_1, J e_1) = 0$.  Consequently, Theorem 3.5 tells us that $\maxmoduli$ consists of the points $C = \bigoplus_{r=1}^t C_r z_r$ $\in \boldgrasstd$ for which the left ideals $C_r \subseteq \la e_1$ satisfy $C_1 \supseteq \cdots \supseteq C_t$;  the ordering of the $C_r$ is due to our choice of Borel subgroup in Section 4, which allowed us to pin down a specific incarnation of $\maxmoduli$. \qed \enddemo

The upcoming example is essentially known.  Variants were presented in \cite{\GeomI, proof of Theorem G} to demonstrate realizability of arbitrary affine varieties as moduli spaces in the classification of uniserial modules; these were then adapted by Hille in \cite{\Hil, Example} to show that any projective variety arises as a moduli space for uniserial modules as well. For the convenience of the reader, we will describe the construction so as to directly target our present goal, rather than give an annotated reference of predecessors with different objectives.
\smallskip

\definition{Examples 5.4}  Every projective variety $V$ is isomorphic to $\maxmoduli$ for some choice of  $\la$, $T$, and $\bd$.  To see this, suppose that $V \subseteq \PP^m$ is defined by $s \ge 0$ nonzero  homogeneous polynomials $h_1, \dots, h_s$ $\in$ $K[X_0, \dots, X_m]$.  To show that $\maxmoduli \cong V$ for suitable $\la$, $T$ and $\bd$, we choose the following quiver $Q$:

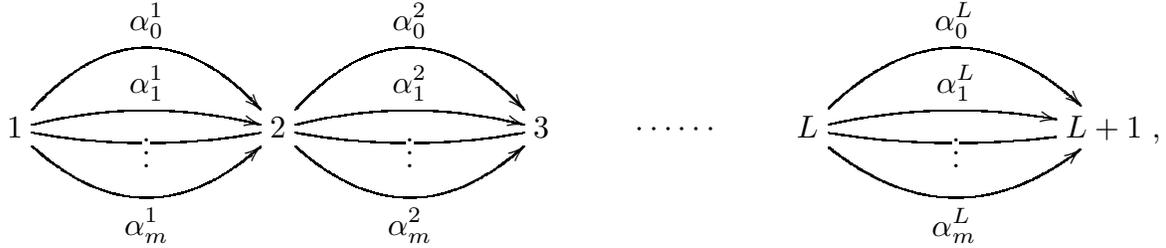
$$\xymatrixrowsep{2.0pc}\xymatrixcolsep{6pc}
\xymatrix{
1 \ar@/^6ex/[r]^{\alpha^1_0} \ar@/^/[r]^{\alpha^1_1} \ar@{}@/_1ex/[r]|{\vdots}
\ar@/_5ex/[r]_{\alpha^1_m}  &2 \ar@/^6ex/[r]^{\alpha^2_0} \ar@/^/[r]^{\alpha^2_1}
\ar@{}@/_1ex/[r]|{\vdots}
\ar@/_5ex/[r]_{\alpha^2_m}  &3 \ar@{}[r]|{\displaystyle\cdots\cdots} &L
\ar@/^6ex/[r]^-{\alpha^{L}_0}
\ar@/^/[r]^{\alpha^{L}_1} \ar@{}@/_1ex/[r]|{\vdots}
\ar@/_5ex/[r]_-{\alpha^{L}_m}  &L+1},$$

\noindent where $L$ is the maximum of the degrees of the $h_r$ if $s \ge 1$,  and an arbitrary positive integer if $s= 0$.  First consider the ideal $I_0 \subseteq  KQ$ which is generated by the relations
$$\alpha_i^{r+1} \alpha_j^r - \alpha_j^{r+1} \alpha_i^r \ \ \ \ \ \text{for} \ \  0 \le i < j \le m\ \ \text{and}\ \  0 \le r < L,$$
and set $\la_0 = KQ/ I_0$.  Observe that, modulo $I_0$, any path $\alpha^v_{i_v} \alpha^{v-1}_{i_{v-1}} \cdots \alpha^u_{i_u}$ in $KQ$  is congruent to the path $\alpha^v_{j_v} \cdots \alpha^u_{j_u}$ with $j_u \le  j_{u+1} \le \cdots \le  j_v$ such that $\{i_u, \dots, i_v\}$ $=$ $\{j_u, \dots, j_v\}$.  To each monomial $X_{i_l} \cdots X_{i_1}$ of degree $l \le L$ in $K[X_0, \dots, X_m]$ we assign the residue class $\alpha^l_{i_l} \cdots \alpha^1_{i_1} + I_0$ in $\la_0$, the order of the variables $X_{i_j}$  being irrelevant by the preceding remark. In this way, each of the polynomials $h_r$ gives rise to an element $f_r + I_0 \in \la_0$ with $f_r \in KQ$.  Letting $I $ be the ideal in $KQ$ which is generated by $I_0$ and $f_1, \dots, f_s$, we define $\la = KQ/I$.  Moreover, we take $T = S_1$ and $\bd = (1, 1, \dots, 1) \in \NN^{L+1}$.  From Example 5.2 we glean that $\maxmoduli = \boldgrasstd$.  In verifying that $\boldgrasstd \cong V$, we will first consider a preliminary case.
\smallskip

\noindent{\bf Case A.}  $s = 0$, i.e., $V = \PP^m$ and $\la = \la_0$.  We will find that, in this case, the map
$$\Psi: \PP^m \rightarrow \boldgrasstd, \ \ \ \ k = [k_0: \dots :k_m] \mapsto C(k):= \sum_{0 \le i < j \le m} \la (k_i \alpha_j^1 - k_j \alpha_i^1)$$
is an isomorphism of varieties.  Given that the assignment $ \PP^m \rightarrow \Gr(m,\sum_{i=0}^m K \alpha_i^1) = \Gr(m,K^{m+1})$ which sends $k = [k_0: \dots :k_m]$ to $\sum_{0 \le i < j \le m} K (k_i \alpha_j^1 - k_j \alpha_i^1)$ is known to be an isomorphism, we will only address well-definedness and surjectivity of $\Psi$.

Suppose $k \in \PP^m$.  In showing that $C(k) \in \boldgrasstd$, it is harmless to assume that $k_0 = 1$, whence $C(k) = \sum_{j=1}^m \la (\alpha_j^1 - k_j \alpha_0^1)$.  Write $M = \la e_1/ C(k) = \la x$, where $x = e_1 + C(k)$.  A straightforward induction on $l \le L$ then yields:
$$(\dagger)   \ \ \  \alpha_{i_l}^l \cdots \alpha_{i_1}^1 x =  k_{i_l} \cdots  k_{i_1} \alpha_0^l \cdots \alpha_0^1 x$$
for any choice of $i_1, \dots, i_l \in \{0, \dots, m\}$.  This shows that $J^l M/ J^{l+1} M \cong K \alpha_0^l \cdots \alpha_0^1 x$, whence we only need to check that $\alpha_0^L \cdots \alpha_0^1 x \ne 0$.  The latter amounts to the easy observation that the path $\alpha_0^L \cdots \alpha_0^1$ does not belong to the left ideal $\sum_{j=1}^m KQ (\alpha_j^1 - k_j \alpha_0^1) + I_0$ of $KQ$.

To verify surjectivity of $\Psi$, let $C \in \boldgrasstd$.  This implies that the intersection $C \cap \sum_{i=0}^m K \alpha_i^1$ has codimension $1$ in $\sum_{i=0}^m K \alpha_i^1$, showing that this intersection is of the form $C(k) \cap \sum_{i=0}^m K \alpha_i^1$ for some $k \in \PP^m$.  One concludes $C = C(k)$.
\smallskip  

\noindent{\bf Case B.}  $s \ge 1$.  Let $\bigl(\boldgrasstd \bigr)_0$ be the variety parametrizing the $\la_0$-modules with top $T$ and dimension vector $\bd$.  Then $\bigl(\boldgrasstd \bigr)_0$ is a copy of $\PP^m$ as we saw in Part A.  Clearly, a point $C = C(k)$ $\in$ $\bigl(\boldgrasstd \bigr)_0$  --  keep  the notation of Part A  --  belongs to $\boldgrasstd$ if and only if the images of $f_1, \dots f_s $ in $\la_0$ belong to $C$.  Therefore our task amounts to showing: $f_r + I_0 \in C(k)$ for all $r \le s$ $\iff$ $k \in V$.  In checking this equivalence, we may again assume that $k_0 = 1$.  In light of equality ($\dagger$) of Part A, we then obtain the congruences $f_r + I_0 \equiv h_r(k) \alpha_0^{\deg f_r} \cdots \alpha_0^1 + I_0$ modulo $C(k)$, whence our claim follows from the fact that we know the image of $\alpha_0^{\deg f_r} \cdots \alpha_0^1$ in $\la_0$ to lie outside $C(k)$.  \qed
 \enddefinition

%%%%%%%%%%%%%%%%%%%%%%%%%%%%%%
%%%%%%%%%%%%%%%%%%%%%%%%%%%%%%

\head 6.  Proof of Theorem 4.4.  Construction of the universal family  \endhead 

We will use Pl\"ucker coordinates for certain affine charts of $\boldgrasstd$, obtained as intersections of suitable open Schubert cells of $\Gr(\dim P - | \bd | ,\, P)$ with $\boldgrasstd$.  In order to introduce convenient cells, we consider ``normalized" linearly independent subsets of $P$ which induce $K$-bases for factor modules of $P$ with dimension vector $\bd$ as follows.  A {\it path\/} in $P$ is any nonzero element $pz_r \in P$, where $p$ is a path in the quiver $Q$ and $r \le t$, i.e., $p$ is a path in $KQ \setminus I$ starting in the vertex $e(r)$ which norms $z_r$; we will say that $p z_r$ ends in the vertex $e_i$ if $p$ does.  In particular, $z_r = e(r) z_r$ is a path in $P$ which ends in $e(r)$.  A {\it path basis with top $T$ and dimension vector $\bd$\/} is any $K$-linearly independent set $\S$ of paths in $P$ such that $z_1, \dots, z_r \in \S$ and precisely $d_i$ of the paths in $\S$ end in $e_i$.  We remark that path bases are generalizations of the skeleta considered in \cite{\generic}.

Clearly, there are only finitely many path bases, and, for any point $C \in \boldgrasstd$, the factor module $P/C$ has a basis $(pz_r + C)_{pz_r \in \S}$ for some path basis $\S$.  The set $\S$ is also referred to as a path basis for $P/C$ in this situation. Given any path basis $\S$ with top $T$ and dimension vector $\bd$, we define  
$$\Schu(\S)  = \{C \in \boldgrasstd \mid P/C \ \text{has basis}\ \S \}.$$
Each of these sets is affine and open in $\boldgrasstd$, because $\Schu(\S)$ is just the intersection of $\boldgrasstd$ with the big Schubert cell $\{C \in \Gr(\dim P - |\bd| ,\, P) \mid P = C \oplus \bigoplus_{pz_r \in \S} K pz_r\}$.  Finally, given any path basis $\S$ and $r \le t$, we denote by $\S_r$ the subset of $\S$ which consists of the paths $pz_r$ in $\S$.  By definition, $z_r \in \S_r$ for each $r$, and $\S$ is the disjoint union of the $\S_r$, $r \le t$.   

The argument for Theorem 4.4 is subdivided into three parts.

\subhead Part A. Claim:  The map $\pi: \maxtopdeg \rightarrow \maxmoduli$ of Section 4 is a morphism of varieties \endsubhead

By Theorem 3.5, the $\autlap$-action boils down to an $\autlat$-action on $\maxtopdeg$.  We again identify $\autlat$ with the usual subgroup of $\autlap$, unless we emphasize identification with $\GL\bigl(\bigoplus_{r=1}^t Kz_r\bigr)$.   Once more, we decompose $P$ in the form $P = \bigoplus_{i=1}^n P_i$, where $P_i = \bigoplus_{r,\, z_r = e_i z_r} \la z_r$.  First, we observe that $\maxtopdeg$ is the disjoint union of the following closed $\autlat$-stable subsets $\frak M (\bd_1, \dots, \bd_n)$, for sequences of dimension vectors $\bd_1, \dots, \bd_n$ adding up to $\bd$: 
$$\frak M (\bd_1, \dots, \bd_n) \ = \ \{C \in \maxtopdeg \mid \underbardim P_i / (C \cap P_i)  = \bd_i \};$$
for justification, recall that all points $C \in \maxtopdeg$ are direct sums of their intersections $C \cap P_i$.  Since the sets $\frak M (\bd_1, \dots, \bd_n)$ are unions of connected components of $\maxtopdeg$, it suffices to show that the restrictions of $\pi$ to these subvarieties are morphisms.  But the restricted maps are direct products of components of the form 
$$\imaxtopdeg \rightarrow \imaxmoduli$$
for $T_i = S_i^{t_i}$, and consequently we do not lose generality in cutting down to the following situation:
\smallskip  

\centerline{{\it Throughout Part A, we assume that $T$ is homogeneous, say  $T =  (S_1)^t$.}}
\smallskip  

\noindent Accordingly,  $P =  P_1 = \bigoplus_{r=1}^t \la z_r$, where each $z_r$ is normed by $e_1$, and our preferred Borel subgroup $\B$ of $\autlat$ (see beginning of Section 4) reduces to the group of upper triangular automorphisms of $T = \bigoplus_{r=1}^t K z_r$, relative to the basis $(z_r)_{r \le t}$.
\smallskip

To prove Claim A, we will use the affine test, applied to the open covers $\maxtopdeg \cap \Schu(\S)$ and $\maxmoduli \cap \Schu(\S)$ of $\maxtopdeg$ and $\maxmoduli$, respectively, where $\S$ traces the $\bd$-dimensional path bases with top $T$.  Since, by Observation 4.1, $\maxtopdeg$ is closed in $\boldgrasstd$, the intersections $\maxtopdeg \cap \Schu(\S)$ indeed constitute an affine cover of $\maxtopdeg$; analogously, the $\maxmoduli \cap \Schu(\S)$ are affine. 

In a first step, we will prune this cover, specializing to more manageable path bases $\S$, so that the corresponding intersections $\maxtopdeg \cap \Schu(\S)$ will still cover $\maxtopdeg$.  We call a $\bd$-dimensional path basis $\S$ with top $T$ {\it totally ordered\/}, if the subsets $\S_r$ of $\S$, for $1 \le r \le t$, are linearly ordered:  By that we mean that the sets $\overline{\S}_r = \{p \text{\ path in\ }  KQ  \mid p z_r \in \S_r\}$ are totally ordered under inclusion.  The $\S_r$ being pairwise disjoint by definition of a path basis, we use the notation ``$\S_r \le \S_s$'' in case $\overline{\S}_r \subseteq \overline{\S}_s$ to avoid confusion.  

The purpose of the following auxiliary observation is to ascertain that the intersections $\maxtopdeg \cap \Schu(\S)$ corresponding to totally ordered path bases $\S$ still cover $\maxtopdeg$.  The argument is somewhat technical, since the open subvarieties $\Schu(\S)$ of $\boldgrasstd$, depending on the distinguished sequence of top elements of $P$, fail to to be stable under the $\autlat$-action in general.

\proclaim{Lemma 6.1}    $\maxtopdeg \ \  \subseteq \ \  \bigcup_{\S \text{\ totally ordered}} \Schu(\S)$. \endproclaim

\demo{Proof}  By Observation 4.3, every $\autlat$-orbit of $\maxtopdeg$ contains a point $C = \bigoplus_{r=1}^t C_r z_r$ with $C_1 \supseteq \cdots \supseteq C_t$.  Clearly, $C$ belongs to $\Schu(\S)$ for some path basis $\S$ satisfying 
$$\S_1 \le \cdots \le \S_t\,.$$
Thus $\S$ is totally ordered.

We proceed to show that every point $D$ in the $\autlat$-orbit of $C$ belongs to $\Schu(\S')$ for some totally ordered path basis $\S'$, which results from $\S$ through a permutation of the sets $\S_r$.  To set up an induction on the number of distinct $C_r$, suppose that there are $u$ distinct left ideals among the $C_r$.  Suppose, moreover, that these are represented by $C_{m_u} \supsetneqq C_{m_{u-1}} \supsetneqq \cdots  \supsetneqq C_{m_1}$, and that each $C_{m_j}$ occurs with multiplicity $s_j$.  Then the subsets $\S_{m_1}, \dots, \S_{m_u}$ of $\S$ are ordered in reverse: $\S_{m_1} > \S_{m_2} > \cdots > \S_{m_u}$.   Finally, for $j \le u$, let $\K_j$ be the subset of $\{1, \dots, t\}$ consisting of the indices $r$ with the property that $C_r = C_{m_j}$; in particular, $|\K_j| = s_j$.
 
Now suppose $D = gC$ for some $g \in \autlat$.  A straightforward induction on $u$ shows that the set $\{1, \dots, t\}$ may be partitioned into subsets $\L_1, \dots, \L_u$ of cardinality $|\L_j| = s_j$, respectively, such that the submodules $L_0 = 0$ and $L_j = \sum_{r \in \bigcup_{s \le j} \L_s}  \la (z_r + D)$ of $P/D$ satisfy the following conditions  for $j = 1, \dots, u$ (we write $\overline{z}$ for residue classes modulo $D$):
\smallskip

$\bullet$ $ L_j  /  L_{j - 1}  = \bigoplus_{r \in \L_j}  \bigl(  \la \overline{z}_r  + L_{j - 1} \bigr)/ L_{j-1}$ for $j \ge 1$;

$\bullet$ $\bigl( \la  \overline{z}_r  + L_{j - 1} \bigr)/ L_{j-1} \ \cong \ \la e_1 / C_{m_j}$ for $r \in \L_j$, with $\overline{z}_r+L_{j-1}$ corresponding to $e_1+C_{m_j}$.
\smallskip

Let $\tau$ be any permutation of $\{1, \dots,  t\}$ which takes $\K_j$ to $\L_j$ for $1 \le j \le u$, and define $\sigma'_r = \sigma_{\tau(r)}$.  Then our choice of $\L_j$ guarantees that the $D$-residue classes of the totally ordered path basis $\S' = \bigcup_{r \le t} \S'_r$ in $P$ constitute a basis for $P/D$. \qed \enddemo

In light of the affine criterion (see, e.g.,  \cite{\Hum, Proposition on p\. 19}), we conclude:

\proclaim{Consequence 6.2}  To prove Claim A, it suffices to show that, for each totally ordered path basis $\S$ with top $T$ and dimension vector $\bd$, the restriction 
$$\pi_{\S}:  \maxtopdeg \cap \Schu(\S) \longrightarrow \maxmoduli \cap \Schu(\S)$$ 
is a morphism of affine varieties. \endproclaim  

The only non-obvious statement of the upcoming lemma is the final assertion under (a), concerning the map $\Psi$.  Expressing suitable restrictions of $\psi$ in terms of Pl\"ucker coordinates is elementary, if somewhat cumbersome; we leave the detail to the reader.  Part (b) of the lemma is obvious, stated solely for easy reference.

\proclaim{Lemma 6.3}  Let $d, k , m , n$ be positive integers with $n \ge k \ge d$ and $n \ge m \ge d$.  Moreover, suppose $V$ is an $n$-dimensional vector space over $K$, $W$ a $k$-dimensional subspace, and $X$ the subset of the classical Grassmannian $\Gr(m, V)$ which consists of the $m$-dimensional subspaces $D$ of $V$ such that $\dim_K (D \cap W) = d$.  Then:
\smallskip

{\rm  {\bf  (a)}}  $X$ is a locally closed subvariety of $\Gr(m,V)$, and the map $\Psi: X \rightarrow \Gr(d, W)$ which sends $D$ to $D \cap W$ is a morphism. 
\smallskip

{\rm  {\bf  (b)}} If $U$ is a complement of $W$ in $V$ and $\pr: V \rightarrow U$ the projection along $W$, then $X$ coincides with the set of all spaces $D \in \Gr(m,V)$ with the property that $\dim_K \pr(D) = m -d$, and the map $\chi: X \rightarrow \Gr(m-d, U)$ which sends $D$ to $\pr(D)$ is a morphism of varieties.
\endproclaim

In the sequel, we will {\it keep a totally ordered path basis $\S$ with top $T$ and dimension vector $\bd$ fixed\/}.  We continue to assume that there are precisely $u$ distinct candidates among the sets $\overline{\S}_r$ consisting of the paths $p$ in $KQ$ such that $pz_r \in \S_r$, respectively;  say
$$\S_{m_1} > \S_{m_2} > \cdots >  \S_{m_u}$$
as before, and denote by  $\L_j$  the set of those indices $r$ in $\{1, \dots, t\}$ for which $\overline{\S}_r = \overline{\S}_{m_j}$.  Moreover, for $j \le u$, we set 
$$Q_j = \bigoplus_{r \in \L_j} \la z_r \,;$$
in particular, $P= \bigoplus_{j=1}^u Q_j$. We write $T_j$ for the top of $Q_j$.  Finally, for any subset $\I \subseteq \{1, \dots, t\}$, we introduce
$$\pr_{r \notin \I}: P \longrightarrow \bigoplus_{r \le t,\, r \notin \I} \la z_r,$$
the canonical projection along the direct summand $\bigoplus_{r \in \I} \la z_r$. 

To show that the map $\pi_\S$ of Consequence 6.2 is a morphism, we will inductively factor it into more transparent components.  In preparation for this factorization, suppose that $C \in \maxtopdeg \cap \Schu(\S)$.  From Theorem 3.5, we know that
$$P/C \cong \bigoplus_{j=1}^u \, \bigoplus_{r \in \L_j}  \la z_r /L_j z_r \ \cong  \ \bigoplus_{j=1}^u \, \bigl(\la e_1/ L_j \bigr)^{\L_j},$$
 where $L_1 \subset \cdots  \subset L_u$ are distinct left ideals contained in $\la e_1$, such that for each $j \le u$ and $r \in \L_j$, the local module $\la z_r/ L_j z_r$ has basis $\S_{m_j}$.  Then $\pi$ sends $C$ to the unique point $C = \bigoplus_{r=1}^t C'_r z_r$ with $C'_1 \supseteq \cdots \supseteq C'_t$ and $P/C \cong P/C'$; in particular, the $L_i$ are the distinct terms on the list $C'_1, \dots, C'_t$, ordered in reverse. 
 
\proclaim{Lemma 6.4}  Let $C \in \maxtopdeg \cap \Schu(\S)$.  Using the notation just introduced,
we obtain:
\smallskip

{\rm {(\bf a)}}  $C \cap Q_1 = \bigoplus_{r \in \L_1} L_1 z_r = \bigoplus_{r \in \L_1} C \cap \la z_r$.  Moreover,  the submodule 
$$C^{(1)} = \bigl( C \cap Q_1\bigr) \oplus \pr_{r \notin \L_1}(C)$$
 of $P$ belongs  to the orbit $\autlat.C$.
 \smallskip
 
{\rm {(\bf b)}} The dimension vector of the intersection $C \cap Q_1$ is $\bold{c}_1 = |\L_1| \cdot (\underbardim  \la e_1 - \underbardim  \S_{m_1})$, and hence does not depend on the choice of $C$. That is, this dimension vector is constant on $\maxtopdeg \cap \Schu(\S)$, as is the dimension vector of $\pr_{r \notin \L_1}(C)$.
\smallskip

{\rm {(\bf c)}}  The map 
$$\psi: \maxtopdeg \cap \Schu(\S) \rightarrow \grass^{T_1}_{\bold{d}_1}, \  \  D \mapsto  D \cap Q_1,$$
where ${\bold{d}_1} = \underbardim Q_1 - \bold{c}_1$, is a morphism of varieties.
\smallskip

{\rm {(\bf d)}}  The map
$$\chi:  \maxtopdeg \cap \Schu(\S) \rightarrow \grass^{\bigoplus_{j \ge 2} T_j}_{\bd - \bold{d} _1},  \ \ D \mapsto \pr_{r \notin \L_1}(D)$$ 
is a morphism.
\endproclaim

\demo{Proof}  {\bf (a)}  By construction, we have an isomorphism $\bigoplus_{j=1}^u \bigoplus_{r \in \L_j} \la z_r / L_j z_r \rightarrow P/C$ which is induced by an automorphism in $\autlat$.   Therefore, each top element  $z_r + C$ of $P/C$ is annihilated by $L_1$, the smallest of the left ideals $L_j$.  This guarantees that   
$$\bigoplus_{r \in \L_1} L_1 z_r\  \subseteq  \ \bigoplus_{r \in \L_1} C \cap \la z_r \ \subseteq \  C \cap Q_1.$$  
Consequently, $Q_1 / (C\cap Q_1)$
has dimension at most $ |\L_1| \cdot (\dim \la e_1 - \dim L_1)$ $=$ $|\L_1| \cdot |\S_{m_1}|$, the latter being the dimension of  $\bigoplus_{r \in \L_1} \bigl( \la z_r /L_1 z_r \bigr)$.   On the other hand, the subset $\bigcup_{r \in \L_1} \S_r \subseteq Q_1$ of cardinality $|\L_1| \times |\S_{m_1}|$ is linearly independent modulo $C$, showing that $Q_1 / (C\cap Q_1)$ has at least this $K$-dimension.  Hence the inclusions displayed above are equalities. 

To see that $P/C^{(1)} \cong P/C$, set $M = P/C$, and note that $P/C^{(1)} \cong U \oplus M/U$, where $U$ is the submodule  of $M$ which is generated by the cosets $(z_r + C)_{r \in \L_1}$. By the preceding paragraph, $U \cong (\la e_1/ L_1)^{\L_1}$. In light of the isomorphism $M\cong \bigoplus_{j=1}^u (\la e_1/L_j)^{\L_j}$, we infer $M \cong U \oplus (\la e_1/ L_2)^{\L_2} \oplus \cdots \oplus (\la e_1/ L_u)^{\L_u}$.  Therefore $M/U \cong  (\la e_1/ L_2)^{\L_2} \oplus \cdots \oplus (\la e_1/ L_u)^{\L_u}$, and $P/C^{(1)} \cong M$ as claimed. 
\smallskip

Part {\bf (b)} is immediate from {\bf (a)}.
\smallskip

Part {\bf (c)} is covered by the first portion of Lemma 6.3, and {\bf (d)} is covered by the second.  \qed
\enddemo 

A straightforward induction will now yield our principal lemma.  In its statement, we retain the notation pertaining to the totally ordered path basis $\S$, as introduced ahead of Lemma 6.4.  In particular, $P = \bigoplus_{j=1}^u Q_j$ with $ Q_j = \bigoplus_{r \in \L_j} \la z_r$.

\proclaim{Lemma 6.5}  Let $C = C^{(0)} \in \maxtopdeg\cap \Schu(\S)$.
\smallskip

For $1 \le j \le u$, we recursively define 
$$C^{(j)} : = \biggl(  C^{(j - 1)} \cap \bigoplus_{i=1}^j Q_i \bigg) \, \oplus\, \pr_{r \notin \L_1 \cup \L_2 \cdots \cup \L_j}(C^{(j-1)}).$$
Then $C^{(j)}$ equals $\bigl( \bigoplus_{r \in L_1 \cup \cdots \cup \L_{j}} (C^{(j-1)} \cap \la z_r) \bigr)\,  \oplus\, \pr_{r \notin \L_1 \cup  \cdots \cup \L_j}(C^{(j-1)})$ and belongs to the $\autlat$-orbit of $C$.  

Moreover, the map $\rho_j: \maxtopdeg \cap \Schu(\S) \rightarrow \maxtopdeg \cap \Schu(\S)$, $C^{(j-1)} \mapsto C^{(j)}$ is a morphism of varieties.
\endproclaim

\demo{Proof by induction on $j \le u$} In case $j = 1$, the first claim is covered by Lemma 6.4; to obtain the final assertion in this case, note that $\rho_1$ is the direct sum of the morphisms $\psi$ and $\chi$ of  that lemma. 
Now let $1 \le j \le u$.  To move from $j-1$ to $j$, we apply Lemma 6.4 to the following pared-down scenario:  Namely, $P$ is replaced by $P^{*} = \bigoplus_{i=j}^u Q_i$, and $T$ by the top $T^{*}$ of  $P^{*}$;  moreover, $C$  is replaced by $C^{*} = \pr_{r \notin \L_1 \cup \cdots \cup \L_{j-1}}(C) =  \pr_{r \notin \L_1 \cup \cdots \cup \L_{j-1}}(C^{(j-1)})$.  Then $ (C^{*} \cap Q_j) \oplus \pr_{r \notin \bigcup_{i \le j} \L_i} (C^{*})$ plays the role of $(C^{*})^{(1)}$, the path basis $\S^{*} = \bigoplus_{r \in \bigcup_{i \ge j} \L_i} \S_r$ of $P^{*}/C^{*}$ takes over the role of $\S$, and the starting conditions of Lemma 6.4 are reproduced for the starred quantities:  Namely, given that $P/C \cong \bigoplus_{i=1}^u (\la e_1/ L_i e_1)^{\L_i}$ with left ideals $L_1 \subsetneqq \cdots \subsetneqq L_u \subsetneqq \la e_1$ has no proper top-stable degenerations, neither does $P^{*}/C^{*}$.   To justify this, we invoke Theorem 3.5 and the fact that  $P^{*}/C^{*} \cong  \bigl(\la e_1/L_{j} \bigr)^{\L_{j}} \oplus \cdots \oplus   \bigl(\la e_1/L_u\bigr)^{\L_u}$ by the induction hypothesis; indeed, the choice of $C$ in $\maxtopdeg$ guarantees that $f(L_i) \subseteq L_k$ for all $f \in \Hom_\la (\la e_1, J e_1)$ and all indices $i,k\in \{1, \dots, t\}$ (cf\. Theorem 3.5 and the subsequent remark).  We apply Lemma 6.4 to this adjusted situation, and re-introduce the summand $C \cap \bigoplus_{i=1}^{j-1} Q_i$ at the end to obtain the claim for $j$.   \qed \enddemo

\proclaim{Corollary 6.6}  The restricted map $\pi_\S: \maxtopdeg \cap \Schu(\S) \rightarrow \maxmoduli$ is a morphism. \endproclaim

\demo{Proof}  Keeping the notation of the preceding lemmata, we let $\tau$ be any permutation of $\{1, \dots, t\}$ such that the elements in $\tau(\L_j)$ precede all of the elements in $\tau(\L_i)$ whenever $i < j$; that is, $L_{\tau(r)} \supsetneqq L_{\tau(s)}$ for all $r \in \L_i$ and $s \in \L_j$ when $i < j$.  Clearly $\tau$ induces an automorphism of $P = \bigoplus_{r=1}^t \la z_r$ that permutes the top elements $z_r$ as prescribed by $\tau$.   Denote by $\overline{\tau}$ the isomorphism $\boldgrasstd \rightarrow \boldgrasstd$ which is induced by this automorphism of $P$.  Then $\pi_\S = \overline{\tau} \circ \rho_u   \circ \cdots \circ \rho_1$.  Indeed, by Lemma 6.4,  the composition $\rho_u   \circ \cdots \circ \rho_1$ leaves the $\autlat$-orbits of $\maxtopdeg \cap \Schu(\S)$ invariant and, for $C \in \maxtopdeg \cap \Schu(\S)$, we have $\rho_u   \circ \cdots \circ \rho_1 (C) = \bigoplus_{r=1}^t C \cap  \la z_r = \bigoplus_{j=1}^u \, \bigoplus_{r \in \L_j} \la z_r/ L_j z_r$.  Applying the map $\overline{\tau}$ to the latter point adjusts the ordering so that $\dim C \cap \la z_r \ge \dim C \cap \la z_s$, whenever $r, s \in \{1, \dots, t\}$ with $r \le s$.  This guarantees that, for any point $C \in \maxtopdeg \cap \Schu(\S)$, the distinguished Borel subgroup $\B$ of $\autlat$ is contained in the stabilizer of $\overline{\tau} \circ  \rho_u   \circ \cdots \circ \rho_1(C)$, meaning that the latter point is the unique candidate $\pi(C)$ in the intersection $\maxmoduli \cap (\autlat.C)$.  Thus $\pi_\S$ factors as claimed, and in light of Lemma 6.5 we conclude that $\pi_\S$ is a morphism.  \qed \enddemo

This completes Part A .  

\subhead Part B. The variety $\maxmoduli$ as a geometric quotient $\maxtopdeg / \autlat$ \endsubhead

It is now easy to prove the ``supplementary information" provided by Theorem 4.4; keep in mind that the $\autlap$-action on $\maxtopdeg$ reduces to the action of $\autlat$.

\proclaim{Proposition 6.7} The map $\pi: \maxtopdeg \rightarrow \maxmoduli$ is a geometric quotient of $\maxtopdeg$ by its $\autlat$-action. \endproclaim

\demo{Proof}  By part A, $\pi$ is a morphism.  That $\pi$ is surjective is obvious from our construction in Section 4, and that its fibers coincide with the $\autlat$-orbits of $\maxtopdeg$ follows from Observation 4.3.  Openness of $\pi$ is an immediate consequence of the facts that $\Im(\pi) = \maxmoduli$ is a subvariety of $\maxtopdeg$ and that the restriction of $\pi$ to this subvariety is the identity:  Indeed given an open subset $W$ of $\maxmoduli$, its closure $\autlat.W = \bigcup_{g \in \autlat}\, g.W$ under the $\autlat$-action is in turn open, and $\pi(W) = \autlat.W  \cap \maxmoduli$.  

Now suppose that $U$ is an open subvariety of $\maxmoduli$.  To check that the comorphism $\pi^0$ induces a bijection from the ring $\Cal O (U)$ of regular functions on $U$ to the ring of those regular functions in $\Cal O \bigl( \pi^{-1} (U) \bigr)$ which are constant on the $\autlat$-orbits of $\pi^{-1} (U)$, let $f \in  \Cal O \bigl( \pi^{-1} (U) \bigr)$ be constant on $\autlat$-orbits.  Then $U = \pi^{-1}(U) \cap \maxmoduli$, and $f = f|_U \circ \pi |_{\pi^{-1}(U)}$.  Uniqueness of this factorization is clear.  \qed 
\enddemo

\subhead Part C.  $\maxmoduli$ as a fine moduli space for our classification problem \endsubhead

We are now in a position to show that $\maxmoduli$ is a fine moduli space, classifying modules of dimension vector $\bd$ which are degeneration-maximal among those with top $T$, up to isomorphism exhibit.  
In light of Proposition 6.7, the construction of a universal family is straightforward.  Indeed, the restriction to $\maxmoduli$ of the following, essentially tautological, family of $\la$-modules indexed by $\boldgrasstd$ will be seen to satisfy our requirements.
\smallskip

\noindent{\bf Construction 6.8 of a family $(\Delta, \delta)$ of modules parametrized by $\boldgrasstd$.}

We again base our construction on the open affine cover $\bigl( \Schu(\S) \bigr)_\S$ of $\boldgrasstd$, where $\S$ runs through the path bases with top $T$ and dimension vector $\bd$.  We will tacitly assume that all the considered $\Schu(\S)$ are nonempty.  For a given path basis $\S$, we let $\Delta_{\sigma} = \Schu(\S) \times K^d$ be the trivial bundle over $\Schu(\S)$, where $d=|\bd|$.  Moreover, we impose and fix an ordering on each such $\S$ --  say $\S$ consists of the elements $b_1 <  \cdots < b_d$ --   and identify the canonical basis of $K^d$ with $\S$ so that the $i$-th canonical basis vector corresponds to $b_i$.  Next we define $\delta_{\S}: \la \rightarrow \End(\Delta_{\S})$ via 
$$\lambda \mapsto  \biggl( (C,\, v) \mapsto (C, \, \rho^{\lambda}_{\S} (C)\cdot v ) \biggr),$$
for $\lambda \in \la$, where $\rho^{\lambda}_{\S} (C)$ is the $d\times d$-matrix that encodes the linear map $P/C \rightarrow P/C$, $x \mapsto \lambda x$, relative to the ordered basis $\{ b + C \mid  b \in \S\}$. One verifies that the map $\delta_{\S}$ is a well-defined $K$-algebra homomorphism:  To guarantee that $\delta_{\S}$ indeed takes values in $\End(\Delta_{\S})$, one checks  that $\rho^{\lambda}_{\S}$ is a morphism $\Schu(\S) \rightarrow \AA^d$ (in doing so, it is convenient to use Pl\"ucker coordinates for $\Schu(\S)$ relative to a basis of $P$ which contains $\S$).  To glue the trivial bundles $\Delta_\S$ together to a bundle over $\boldgrasstd$, we use the following morphisms $g_{\S, \S'}: \Schu(\S) \cap \Schu(\S') \rightarrow \GL_d$: Namely,  $g_{\S, \S'}$ sends any $C$ in the intersection to the transition matrix from the ordered basis $\{b + C \mid b \in \S \}$ of $P/C$ to the ordered basis $\{b' + C \mid b' \in \S' \}$ of $P/C'$. Clearly, the maps $g_{\S, \S'}$ then satisfy the relevant cocycle condition.  The resulting vector bundle $\Delta$ over $\boldgrasstd$ carries a $\la$-structure by way of the algebra homomorphism $\delta: \la \rightarrow \End(\Delta)$ which is induced by the $\delta_\S$, the latter maps being compatible with the gluing of the $\Delta_\S$.  
\smallskip    

For the final link in the proof of Theorem 4.4, we now merely have to pull the established facts together. We follow Newstead's blueprint in doing so.

\proclaim{Proposition 6.9} $\maxmoduli$ is a fine moduli space for the modules of dimension vector $\bd$ which are degeneration-maximal among those with top $T$.  The corresponding universal family is the restriction of  the family $(\Delta, \delta)$ constructed in 6.8 to the subvariety $\maxmoduli$ of $\boldgrasstd$.  \endproclaim   

\demo{Proof}  Proposition 6.7 makes $\pi: \maxtopdeg \rightarrow \maxmoduli$ an orbit map which, moreover, is a categorical quotient of $\maxtopdeg$ by the $\autlat$-action; the latter means that $\pi$ has the universal property  stated ahead of Corollary 4.5.  Criterion 2.6 in \cite{\class}, which amounts to a module-theoretic reformulation of Newstead's \cite{\New, Proposition 2.13}, thus guarantees that $\maxmoduli$ is a coarse moduli space for our problem.  

We will apply  \cite{\New, Proposition 1.8} to verify that $\maxmoduli$ is even a fine moduli space.  Let $(\Gamma, \gamma)$ be the restriction of the above family $(\Delta, \delta)$.  Evidently, $(\Gamma, \gamma)$ satisfies Condition (i) of \cite{\New, Proposition 1.8}.  That Condition (ii) is met as well is immediate from our definition of equivalence of families of modules parametrized by a variety; see Section 1. \qed \enddemo

This completes the proof of Theorem 4.4.

\enddocument